\RequirePackage{fix-cm}
\documentclass[smallextended]{svjour3}       
\smartqed  
\usepackage{graphicx}
\usepackage{epsfig}
\usepackage{amsmath}
\usepackage{amssymb}
\usepackage{amsfonts}
\usepackage{mathabx}
\usepackage{euscript}
\usepackage{textcomp}
\usepackage{tikz}
\usetikzlibrary{shapes,arrows}
\usetikzlibrary{shadows,calc}
\usepackage{xcolor}

\usepackage[utf8]{inputenc}
\usepackage{epstopdf}
\usepackage{multirow}
\usepackage{subfig}

\newcommand{\R}{\mathbb{R}}
\newcommand{\bea}{\begin{eqnarray}}
\newcommand{\eea}{\end{eqnarray}}
\newcommand{\mbb}[1]{{\mathbb{#1}}}
\newcommand{\bfm}[1]{{\mbox{\boldmath{$#1$}}}}

\begin{document}

\title{Hybrid BSQI-WENO Based Numerical Scheme for Hyperbolic Conservation Laws
}


\author{Rakesh Kumar         \and
        S. Baskar 
}


\institute{ Rakesh Kumar \at
              Indian Institute of Technology Bombay, Powai, Mumbai 400 076. INDIA\\
\email{rakeshk@math.iitb.ac.in}           
           \and
           S. Baskar \at
              Indian Institute of Technology Bombay, Powai, Mumbai 400 076. INDIA\\
              \email{baskar@math.iitb.ac.in}   }

\date{Received: date / Accepted: date}

\maketitle

\begin{abstract}
In this paper, we intend to use a B-spline quasi-interpolation (BSQI) technique to develop higher order hybrid schemes for conservation laws.  As a first step, we develop cubic and quintic B-spline quasi-interpolation based numerical methods for hyperbolic conservation laws in 1 space dimension, and show that they achieve the rate of convergence 4 and 6, respectively.  Although the BSQI schemes that we develop are shown to be stable, they produce spurious oscillations in the vicinity of shocks, as expected.  In order to suppress  the oscillations, we conjugate the BSQI schemes with the fifth order weighted essentially non-oscillatory (WENO5) scheme.  We use a weak local truncation based estimate to detect the high gradient regions of the numerical solution.  We use this information to capture shocks using WENO scheme, whereas the BSQI based scheme is used in the smooth regions. For the time discretization, we consider a strong stability preserving (SSP) Runge-Kutta method of order three. At the end, we demonstrate the accuracy and the efficiency of the proposed schemes over the WENO5 scheme through numerical experiments.
\end{abstract}
\section{Introduction}
\label{intro}

In this work, we intend to use B-Spline Quasi-Interpolation (BSQI) to develop numerical schemes for the hyperbolic conservation laws  
\begin{equation}\label{hyp.1d}
 \bfm{u}_t+\bfm{f}(\bfm{u})_x=0, \;\;\;(x,t)\in (a,b) \times (0,T],
\end{equation}
where $\bfm{u}:(a,b)\times (0,T]\rightarrow \mbb{R}^p$ and $\bfm{f}:\mbb{R}^p\rightarrow \mbb{R}^p$ is the flux function, along with the initial condition
\begin{equation}\label{hyp1d.ini}
 \bfm{u}(x,0)=\bfm{u}_0(x),\;\;\;x\in [a,b].
\end{equation}
It is well-known that the classical solution of \eqref{hyp.1d} may ceases to exist in finite time, even the initial data is sufficiently smooth. The higher order finite difference, spectral, and finite element methods are well suitable in approximating the solution in the smooth regions, but fail to provide a non-oscillatory solution in the presence of jump discontinuities. Thus, it is challenging to devise an higher order numerical scheme, which preserves the higher order accuracy in smooth regions and produce non-oscillatory solution in a vicinity of shocks. 

The BSQI approximations are known (see Sablonni\`{e}re \cite{sab_05a}) to be very accurate in approximating smooth functions and their derivatives. This motivated many researchers to use this technique to develop numerical  methods for some partial differential equations.  For instance, Li {\it et al.} \cite{li-etal_14b} developed a  mesh-free method for the viscous Burgers' equation, and Zhu and Kang \cite{zhu-kan_10a} developed a finite difference  method for the Burgers-Fisher equation. Recently, we (Kumar and Baskar \cite{kum-bas_16a}) have developed higher order schemes for some Sobolev type equations with special reference to equal width and BBMB equations.  But these equations involve dissipative and dispersive effects that lead to smeared shocks and dispersive shocks, respectively.  Zhu and Kang \cite{zhu-kan_10b} made an initial effort in devising cubic B-spline quasi-interpolation (CBSQI) based numerical scheme for hyperbolic conservation laws.  But their study is restricted only to scalar conservation laws and also they used quasi-linear form of the equation to develop the numerical scheme and presented numerical results for Burgers' equation up to the critical time at which the shock forms. To our knowledge, there is no work on developing a BSQI based conservative numerical method for hyperbolic conservation laws.

In the present work, our interest is to develop higher order methods for \eqref{hyp.1d}-\eqref{hyp1d.ini} using cubic and quintic BSQI (denoted by CBSQI and QnBSQI, respectively) for approximating space derivative and the 3$^{\rm rd}$ order Runge-Kutta method for time derivative. Although the resulting schemes are conditionally stable, as they are higher order methods, they develop spurious oscillations near the shocks.  It is known that such a problem is common for any higher order methods. There are at least two approaches to minimize the non-physical oscillations in the numerical solution, namely, using limiters and using essentially non-oscillatory (ENO) reconstructions.  Out of these, the ENO approach can minimize the oscillations without much compromise in the order of accuracy.  

An improved version of the ENO approach is the well known Weighted ENO (WENO) reconstruction introduced by Liu {\it et al.} \cite{liu-etal_94a} (also see Jiang and Shu \cite{jia-shu_96a}). However, this reconstruction algorithm is expensive. In order to minimize the expense, there are many hybrid schemes introduced in the literature for the hyperbolic conservation laws. For instance, a non-conservative hybrid method is introduced by Harabetian and Pego \cite{har-peg_93a}, the multi-domain methods in one and two space dimensions are developed using Fourier continuation and WENO scheme by Shahbazi {\it et al.} \cite{sha-etal_11a,sha-etal_13a}. Costa and Don \cite{cos-don_07a,cos-don_07b} obtained a high order hybrid method by combining the spectral and the central schemes with WENO.   Cheng and Liu \cite{che-etal_13a,che-liu_14a} introduced multi-domain RKDG methods with proper boundary treatments. The characteristics hybrid compact-WENO scheme is proposed by Ren {\it et al.} \cite{ren-etal_03a}.  

In developing efficient hybrid schemes, choosing smooth indicator is an important component of the algorithm. Karni {\it et al.} \cite{kar-etal_02a}, used a Weak Local Truncation Error (WLTE) based smooth indicator in developing an adaptive algorithm for hyperbolic systems. Further Kurganov and Liu \cite{kur-liu_12a} developed adaptive artificial viscosity methods for the system of Euler equations of motion using Weak Local Residual (WLR) (similar to WLTE).  The WLR is used in many other problems, for instance, Mungkasi {\it et al.} \cite{mun-etal_14a} used WLR for the system of shallow water equations  and Chen {\it et al.} \cite{che-etal_13b} developed the adaptive artificial viscosity method for the {S}aint-{V}enant system using this indicator. To reduce the computational cost in the case of system of Euler equations,  Dewar {\it et al.} \cite{dew-etal_15a} introduced WLR pressure based smooth indicator. A comparison study of many switch indicators for hyperbolic conservation laws is performed by Li and Qiu \cite{li-qiu_10a,li-qiu_14a} for uniform and  curvilinear meshes.

In this paper, we propose higher order hybrid BSQI-WENO schemes for hyperbolic conservation laws.
 In a BSQI-WENO scheme, we conjugate the BSQI based numerical scheme with a WENO scheme. To detect the high gradient region
of the numerical solution, we use a WLR based estimate of Kurganov and Liu \cite{kur-liu_12a}.  We use this information to capture shocks using WENO scheme, whereas the BSQI based scheme is
used in the smooth region. For the time discretization, we consider a SSP Runge-Kutta method of 
order three. To validate the numerical scheme, we perform the numerical experiments using Burgers and Buckley-Leverett equations in the case of scalar conservation laws and Euler equations in the case of system of equations.  It is also interesting and important to study the efficiency and the accuracy of the hybrid methods constructed using the BSQI based numerical schemes and WENO schemes. 

The outline of the paper is as follows. In section \ref{BSQI.BASIC}, we discuss the BSQI and its application in approximating the derivative of a function. In section \ref{bsqi1}, we develop the semi-discrete form of  hyperbolic conservation laws using CBSQI and QnBSQI approximations to the space derivative. We further perform the von Neumann stability analysis and provide the numerical evidence for the order of accuracy.  We develop two hybrid 
 numerical schemes, namely, the Hybrid4 and Hybrid6 by combining the CBSQI with WENO3 and QnBSQI with WENO5, respectively, in section \ref{hybrid.sec}. We perform some numerical experiments in section \ref{numericalhybrid}, to validate the schemes for scalar as well as system of equations. Finally, we test the accuracy and the efficiency of the proposed algorithms over the pure WENO in section \ref{ordeff.sec}.

\section{B-Spline Quasi-Interpolation}\label{BSQI.BASIC}

In this section, we briefly discuss the construction of BSQI of degree upto 5 and give the numerical rate of convergence in each cases.  For more details on the basic concepts of spline spaces and B-splines we refer to Schumaker \cite{sch_07a}, De boor \cite{deb_01a}, N{\"u}rnberger \cite{nur_89a}, and Kvasov \cite{kva_00a}.

Let the space interval $I=[a,b]$ be divide into $m$ subintervals $I_j=[x_{j-1},x_j]$, $j=1,2,\ldots,m$ of equal length $\Delta x=x_j-x_{j-1}$. The boundary points of the subintervals are called the {\it knots} and let $$X_m=\{x_j: j=0,1,\cdots, m\}$$ denotes the set of knots. The {\it spline space} of degree $d$ is defined by 
(see Schumaker \cite{sch_07a})
\begin{equation}
 S^d(I, X_m)=\{ v\in C^{d-1}: v|_{[x_{j-1},x_j]}\in \mbb{P}_d, ~1\leq j\leq m\},
\end{equation}
where $\mbb{P}_d$ is the polynomial space of degree $d$. The set $S^d(I,X_m)$ forms a vector space of dimension $m+d$ and the set  
$
 \{B_j^d: 1\leq j\leq m+d\}
$
of B-splines of degree $d$ forms a basis for this space. B-splines can be constructed recursively by starting with the B-spline of degree 0 given by
\begin{equation}\label{roe}
B_j^0=\left\{\begin{array}{ll}
                1 ~~~~  x\in[x_{j-1},x_{j}], \\
                0~~~~~~ \text{otherwise}.
                \end{array}\right.
\end{equation}
The recursive formula for the higher degree B-splines takes the form
\begin{equation}
{B}_j^d(\xi)=\Big(\frac{x_j-\xi}{x_j-x_{j-d}}\Big){B}_j^{d-1}(\xi)+\Big(\frac{\xi-x_{j-d-1}}{x_{j-1}-x_{j-d-1}}\Big){B}_{j-1}^{d-1}(\xi),
\;\;\xi\in \mathbb{R}.
\end{equation}
The general form of the B-spline quasi-interpolant (BSQI) of degree $d$ for a function $u:I\rightarrow \R$ is given by  (see Sablonni\`{e}re \cite{sab_05a})
\begin{equation}\label{QI.form}
 Q_du(x)=\displaystyle\sum\limits_{j=1}^{m+d}\mu_j(u) {B_j^d}(x),
 \end{equation}
where the coefficients $\mu_j$ can be obtained by imposing the condition that the quasi-interpolant $Q_du$ is exact if $u$ is a polynomial of degree less than or equal to $d$. The coefficients $\mu_j$ are the linear combination of function (approximate) values $u_j\approx u(\xi_j)$ at the nodes, where
\begin{enumerate}  
\item if $d$ is odd, then the nodes are simply the knots, {\it i.e.} $\xi_j=x_j$, $j=0,1,\ldots,m$ and the set of nodes involved in calculating $\mu_j$ is 
\begin{equation}
  \mathcal{N}_j^d:=\{x_{j-d},x_{j-d+1},\ldots, x_{j-1}\},\;\;j= 0,1,\ldots,m+d
 \end{equation}
 and
 \item if $d$ is even, then the nodes are taken as $\xi_j = {(x_{j-1}+x_{j})}/{2}$ and the set of nodes involved in calculating $\mu_j$ is 
 \begin{equation}\label{even.1}
 \mathcal T_j^d=\{\xi_l:\;\;\xi_l={(x_{l-1}+x_{l})}/{2},\;\;x_l\in\mathcal{N}_j^d\cup\{x_{j-d-1},x_j\}\},\;\;j= 1,\ldots,m+d.
\end{equation}
\end{enumerate}
For more details, we refer to Sablonni\`{e}re \cite{sab_05a}.

In case when support of B-spline 
 contains distinct knots, the coefficients $\mu_j$ for quadratic (QBSQI), cubic (CBSQI), quadric (QdBSQI), and quintic (QnBSQI) are given as follows:\\
 \noindent\textbf{Quadratic:}
\begin{equation}
 \mu_j(u)=\frac{1}{8}(-u_{j-1}+10u_j-u_{j+1}) ~~~\mbox{for }~~j=3,\ldots,m,
\end{equation}
\textbf{Cubic:}
\begin{equation}
 \mu_j(u)=\frac{1}{6}(-u_{j-3}+8u_{j-2}-u_{j-1}) ~~~\mbox{for }~~j=3,\ldots,m+1,
\end{equation}
\textbf{Quadric:}
\begin{equation}
 \mu_j(u)=\frac{47}{1152}(u_{j-4}+u_j)-\frac{107}{288}(u_{j-3}+u_{j-1})+\frac{319}{192}u_{j-2}~~~\mbox{for }~~j=5,\ldots,m,
\end{equation}
\textbf{Quintic:}
\begin{equation}
 \mu_j(u)=\frac{13}{240}(u_{j-5}+u_{j-1})-\frac{7}{15}(u_{j-4}+u_{j-2})+\frac{73}{40}u_{j-3}~~~\mbox{for }~~j=5,\ldots,m.
\end{equation}
We have the following theorem for general BSQI in approximating a function:
\begin{theorem}\label{thm.cubic}
There exists a constant $0<C<1$ such that for all  $u\in W^{d+1,\infty}(I)$ and for any set of knots $\xi_m$ on $I$ with step size 
$\Delta x$, we have
\bea
\|u-Q_du\|_\infty \le C\Delta x^{d+1}\|u^{(d+1)}\|_\infty.\nonumber
\eea
\end{theorem}
For the proof of this theorem, we refer to Sablonni\`{e}re \cite{sab_07a}, and DeVore and Lorentz \cite{dev-lor_93a}.

From the construction we can see that a B-splines of degree $d$ is $d-1$ times continuously differentiable.  Thus, we can also approximate the derivatives of a function $u$ with the corresponding
derivatives of BSQI as long as it is possible, {\it i.e.}
\begin{equation}\label{QI.der}
 u^{(l)}(x)\approx (Q_du)^{(l)}(x)=\displaystyle{\sum_{k=1}^{m+d} \mu_k (u)\Big(B_k^d(x)\Big)^{(l)}},~~~~~l\in\{1, 2,\ldots, d-1 \}.
\end{equation}
\begin{table}[t]
\begin{center}
\caption{The coefficient $b_{\cdot,l}$ for $l=0,1$ and for interior nodes, in the case of quadratic, cubic, quadric, and quintic B-spline quasi-interpolations.
}
\begin{tabular}{|c|c| c|c|c|c|c|c|c|c|c|}
\hline
Approximation&B-Spline&$b_{j-4,l}$&$b_{j-3,l}$ & $b_{j-2,l}$
&$b_{j-1,l}$ &$b_{j,l}$ &$b_{j+1,l}$ &$b_{j+2,l}$ &$b_{j+3,l}$
&$b_{j+4,l}$ \\ [0.3ex]

\hline
Function&Quadratic&0 &0 &$\frac{-1}{64}$ &$\frac{1}{16}$
&$\frac{29}{32}$ &$\frac{1}{16}$ &$\frac{-1}{64}$ & 0&0 \\
\cline{2-11}
 ($l=0$) &Cubic&0 &0 &$\frac{-1}{36}$ &$\frac{1}{9}$ &$\frac{15}{18}$
&$\frac{1}{9}$ & $\frac{-1}{36}$& 0&0 \\\cline{2-11}
&Quadric&$\frac{47}{442368}$ &$\frac{131}{18432}$ &
$\frac{-4951}{110592}$&$\frac{6271}{55296}$ &$\frac{62543}{73728}$
&$\frac{6271}{55296}$ &$\frac{-4951}{110592}$
&$\frac{131}{18432}$ &$\frac{47}{442368}$ \\\cline{2-11}
&Quintic &$\frac{13}{28800}$ &$\frac{113}{14400}$ &
$\frac{-101}{1800}$&$\frac{2111}{14400}$ &
$\frac{2311}{2880}$&$\frac{2111}{14400}$ &$\frac{-101}{1800}$
&$\frac{113}{14400}$ &$\frac{13}{2880}$ \\ \hline

 1$^{\rm st}$ derivative &Quadratic &0 &0
&$\frac{1}{16}$&$\frac{-5}{8}$ &0 &$\frac{5}{8}$ &$\frac{-1}{16}$ &0
&0
\\ \cline{2-11}
($l=1$) &Cubic&0 &0 &$\frac{1}{12}$ &$\frac{-2}{3}$ &0 &$\frac{2}{3}$
&$\frac{-1}{12}$ &0 &0 \\\cline{2-11}
&Quadric &$\frac{-47}{55296}$ &
$\frac{-101}{9216}$&$\frac{3751}{27648}$ &$\frac{-20323}{27648}$ &0 &
$\frac{20323}{27648}$&$\frac{-3751}{27648}$
&$\frac{101}{9216}$& $\frac{47}{55296}$ \\ \cline{2-11}
 &Quintic &$\frac{-13}{5760 }$ &
$\frac{-1}{320}$&$\frac{341}{2880}$&$\frac{-2069}{2880}$
&0&$\frac{2069}{2800}$ &$\frac{-341}{2800}$&$\frac{1}{320}$&
$\frac{13}{5760}$\\ \hline
\end{tabular}
\label{BSQI.DER}
\end{center}
\end{table}
The approximation of function $u$ and its first derivative at point $x_j$  are as follows:\\
\noindent \textbf{Case I}. When $d$ is even, we have
\begin{equation}\label{odd}
 (Qu)^{(l)}(x_j)=\frac{1}{\Delta x^l}\displaystyle{\sum_{k=j-d}^{j+d}b_{k,l} u_k},~~~~~~l=0, \ldots,d-1,
\end{equation}
and\\
\textbf{Case II}. when $d$ is odd, we have
\begin{equation}\label{even}
 (Qu)^{(l)}(x_j)=\frac{1}{\Delta x^l}\displaystyle{\sum_{k=j-(d-1)}^{j+(d-1)}b_{k,l} u_k},~~~~~~l=0, \dots,d-1,
\end{equation}
where coefficients $b_{k,l}$ are to be obtained from \eqref{QI.der} and the derivatives of the B-splines.  For $l=0,1$ these coefficients are listed in Table \ref{BSQI.DER}.

\begin{figure}
 \includegraphics[width=8.3cm]{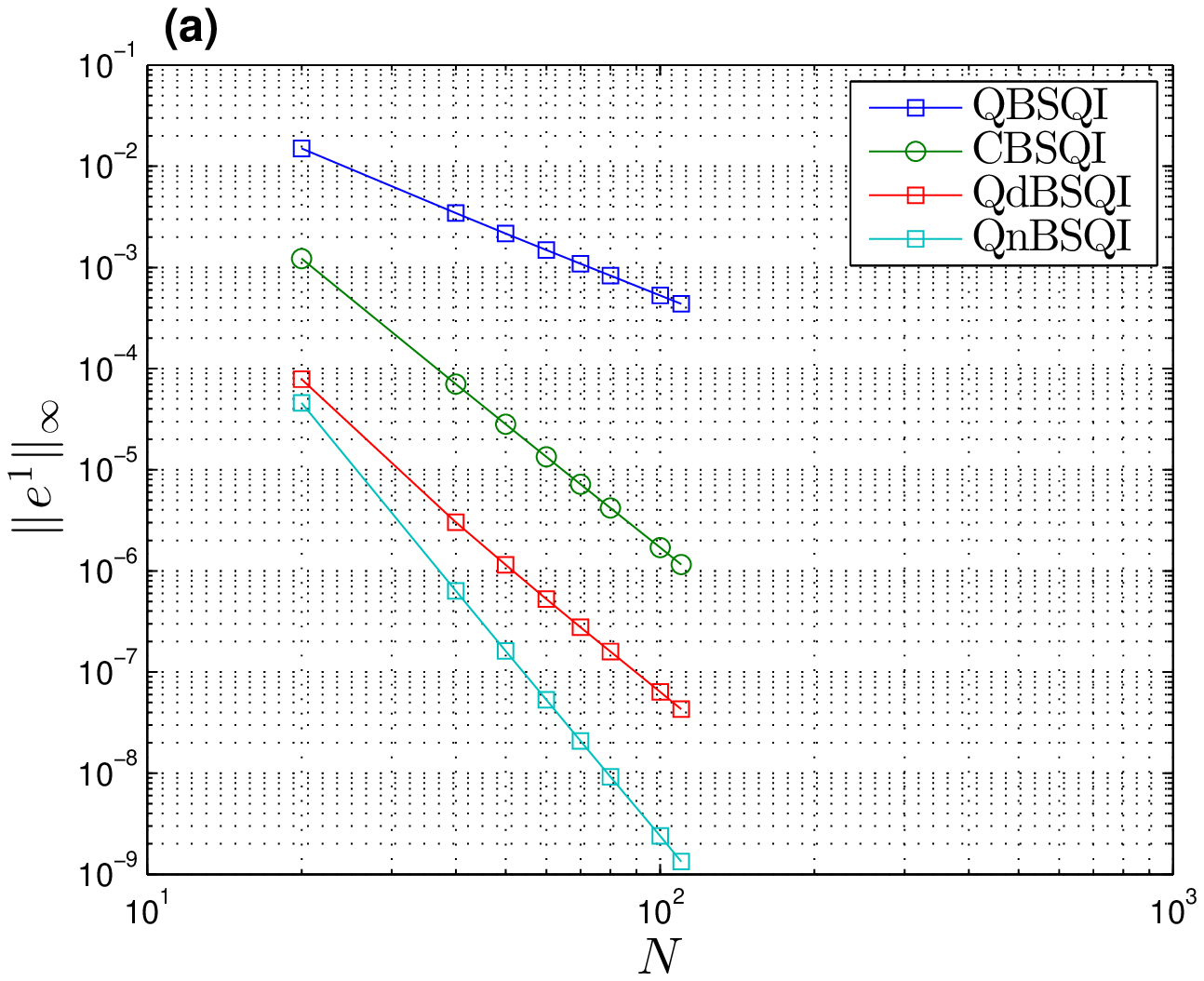}
 \includegraphics[width=8.3cm]{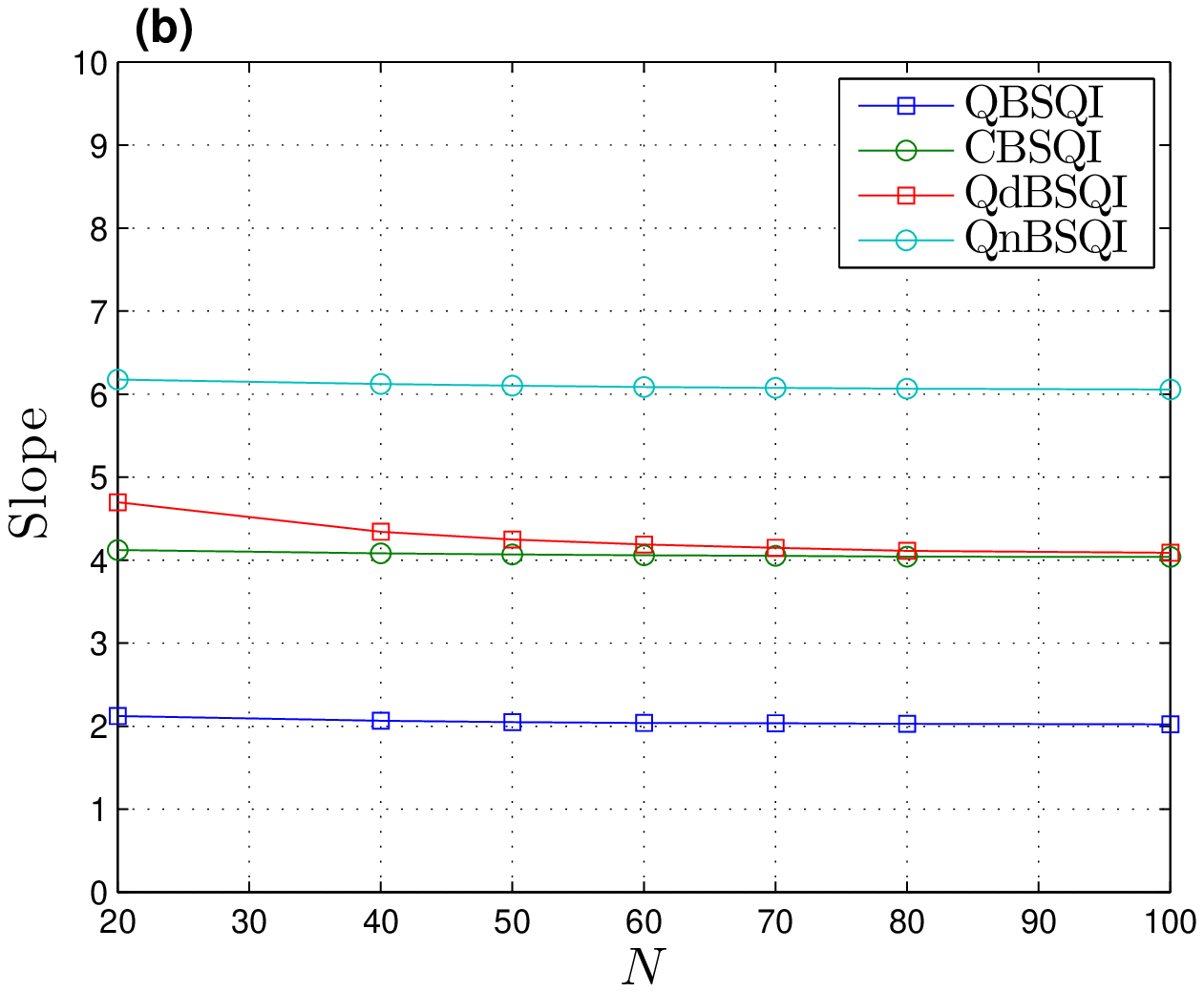}
 \caption{Comparison of BSQI of degree from quadratic to quintic for the approximation of first derivative for model problem \eqref{exm1} in term of $L^{\infty}$-error,
 (a) log-log plot of $L^{\infty}$-error with number of mesh points, (b) Slope of the corresponding BSQI.}
 \label{fir_der}
\end{figure}

In the following example, we demonstrate the accuracy and the rate of convergence in 
  approximating the first derivative of a smooth function.
\begin{example}\label{ex.ordac}{\rm

Consider the function
\begin{equation}\label{exm1}
 u(x)=\sin(\pi x), ~x\in [-1,1],
\end{equation}
whose first derivative is given by
$
 u'(x)= \pi \cos(\pi x).\;\;
$
In order to make the comparison among the BSQI's in terms of their accuracy, we consider the $L^{\infty}$-error norm given by
\begin{equation}
 \| e^k\|_{\infty}=\displaystyle{\max_{j} |u^{(k)}_j-(Q_du)^{(k)}(x_j)|},  \;\; k=0,1,
\end{equation}
where the superscript $(k)$ denotes the $k^{th}$ derivative and $e^k$ denotes the error associated with the $k^{th}$ derivative.
In Figure \ref{fir_der}(a), we depict the log-log plot of the $L^{\infty}$-error and the number of mesh points, in the case of the first derivative approximation of the model problem \eqref{exm1}. The QBSQI converges to the first derivative with  rate two, whereas CBSQI and QdBSQI achieves the fourth order of convergence, which is clear from Figure \ref{fir_der}(b). The QnBSQI has a higher order accuracy among the considered BSQI, which has the rate of convergence six in the case of the first derivative approximation.
}
\end{example}

\section{BSQI Scheme}\label{bsqi1}

In this section, we discuss the semi-discrete form of \eqref{hyp.1d} obtained using the BSQI based spatial approximation.  For the sake of notation simplicity, we restrict to $p=1$, and consider the scalar conservation laws of the form 
\begin{equation}\label{scl.1d}
{u}_t+{f}({u})_x=0, \;\;\;(x,t)\in (a,b) \times (0,T].
\end{equation}
The BSQI formulation for a general 1D system \eqref{hyp.1d} can be done componentwise in a similar way. 

We treat the flux function ${f}$ in \eqref{scl.1d} as $f({u})=f(\xi,t)$ and use the BSQI approximation \eqref{QI.form} in \eqref{scl.1d} to get 
\begin{equation}\label{semidisc.1}
 u_t+\displaystyle{\sum_{k=1}^{m+d}\mu_k(f)(B^d_k)'}\approx 0,\;\;\;x\in [a,b],~t\in (0,T].
\end{equation}
Let us use the notation $u_j\approx u(x_j,t)$, for $x_j\in X_m$, to denote the approximation of the solution $u$ obtained using \eqref{semidisc.1} and write the semi-discrete form of \eqref{scl.1d} as
\begin{equation}\label{semidisc.2}
 \frac{du_j}{dt}=-\displaystyle{\sum_{k=1}^{m+d}\mu_k(f)\Big(B^d_k(x_j)\Big)'}, ~j=d-1,\ldots,m-d+1.
\end{equation}
The above methods can work for any degree of BSQI.  But in this work, we intend to use only the CBSQI and QnBSQI formulations, which are given, respectively, by
\begin{equation}\label{semidisc.CBSQI}
 \frac{du_j}{dt}=-\frac{1}{\Delta x}\Big(\frac{1}{12}f_{j-2}-\frac{2}{3}f_{j-1}+\frac{2}{3}f_{j+1}-\frac{1}{12}f_{j-2}\Big),
\end{equation}
for $j=2,\ldots,m-2$ and
\begin{multline}\label{semidisc.QnBSQI}
 \frac{du_j}{dt}=-\frac{1}{\Delta x}\left(-\frac{13}{5760}f_{j-4}-\frac{1}{320}f_{j-3}+\frac{341}{2880}
 f_{j-2}^n-\frac{2069}{2880}f_{j-1}\right.\\
\left. +\frac{2069}{2880}f_{j+1}-\frac{341}{2880}f_{j+2}+\frac{1}{320} f_{j+3}+\frac{13}{5760}f_{j+4}\right),
\end{multline}
for $j=4,\ldots,m-4$, where we have used the notation $f_j = f(u_j,t)$.
The above semi-discrete formulations can be written in the conservation form as
\begin{equation}\label{semidisc.consform}
 \frac{du_j}{dt}=-\frac{1}{\Delta x}\left(F_{j+\frac{1}{2}} - F_{j-\frac{1}{2}}\right)=:L(u),
 \end{equation}
 where $F_{j\pm\frac{1}{2}}=F_{j\pm\frac{1}{2}}(t)$ denotes the numerical flux given by
 \begin{equation} \label{numflux.CBSQI}
F_{j+\frac{1}{2}}=-\frac{1}{12}f_{j-1}+\frac{7}{12}f_{j}+\frac{7}{12}f_{j+1}-\frac{1}{12}f_{j+2}=: {F}^{\rm CBSQI}_{j+\frac{1}{2}}(t),
\end{equation}
for the CBSQI scheme \eqref{semidisc.CBSQI} and
\begin{multline} \label{numflux.QnBSQI}
 F_{j+\frac{1}{2}}=\frac{1}{\Delta x}\left(\frac{13}{5760}f_{j-3}+\frac{31}{5760}f_{j-2}-\frac{651}{5760}f_{j-1}+\frac{3487}{5760}f_{j}\right.\\
\left. +\frac{3487}{5760}f_{j+1}-\frac{651}{5760}f_{j+2}+\frac{31}{5760}f_{j+3}+\frac{13}{5760}f_{j+4}\right)=:{F}^{\rm QnBSQI}_{j+\frac{1}{2}}(t),
 \end{multline}
 for the QnBSQI scheme \eqref{semidisc.QnBSQI}. 
 
Observe that \eqref{semidisc.consform} forms a system of ODEs.  When the flux $f$ is linear, then \eqref{semidisc.consform} leads to a system of linear ODEs.  In this case, by assuming that the solution of \eqref{semidisc.consform} in the form
\begin{equation}\label{von.sol}
 u_j(t)=\hat{u}(t) e^{ij\theta}, ~j=d-1,\ldots, m-d+1,
\end{equation}
where $i=\sqrt{-1}$ and $\theta=\omega \Delta x$ with $\omega$ being the wave-number, we get
\begin{equation*}
 \frac{du_j}{dt}=-\frac{1}{\Delta x}(\mathcal{C}(\theta) u_j), ~j=d-1,\ldots, m-d+1,
\end{equation*}
where $\mathcal{C}(\theta)$ is given by
\begin{equation*}
 \mathcal{C}(\theta)=
 \left\{\begin{array}{ll}
 \dfrac{i}{6}\big(8\sin(\theta)-\sin(2\theta)\big)&\mbox{for CBSQI},\medskip\\
 \dfrac{i}{6}\Big(\dfrac{2069}{240}\sin(\theta)-\dfrac{341}{240}\sin(2\theta)+\dfrac{3}{80}\sin(3\theta)+\dfrac{13}{480}\sin(4\theta)\Big)&\mbox{for QnBSQI}.
 \end{array}\right.
\end{equation*}

For the time discretization, we use the Strong Stability Preserving (SSP) Runge-Kutta method of order three  given by (see \cite{got-etal_01,rur-spi_04} for more details)
 \begin{align*}
  v^{(0)}&=u^n,\\
  v^{(1)}& =v^{(0)}+\Delta t L(v^{(0)}),\\
  v^{(2)}&=\frac{3}{4}v^{(0)}+\frac{1}{4}(v^{(1)}+\Delta t L(v^{(1)})),\\
  v^{(3)}&=\frac{1}{3}v^{(0)}+\frac{2}{3}(v^{(2)}+\Delta t L(v^{(2)})),\\
  u^{n+1}&=v^{(3)}.
 \end{align*}
  Note that the spectrum of the Runge-Kutta method of order three is known to contain some part of the imaginary axis.  Thus, this method along with the CBSQI and QnBSQI space discretizations can be stable under the condition 
\begin{equation*}
 -\beta\leq \lambda \mathcal{C}(\theta) \leq \beta,
\end{equation*}
where $\beta= \sqrt{3}$ and $\lambda=\Delta t/\Delta x$ (see Hirsch \cite{hir_07a} for more details).  This implies that, we can find the condition on $\lambda$ such that the scaled spectrum of CBSQI and QnBSQI schemes are contained in the stability spectrum of the Runge-Kutta method of order three. 

In the following example, we demonstrate the order of accuracy of the CBSQI and QnBSQI schemes in the case of linear advection equation.
\begin{example}\label{ex.ordacc.adv}{\rm

Consider the linear advection equation 
\begin{equation}\label{adv.11}
 u_t+u_x=0\;\;\; (0,2\pi)\times (0,T],
\end{equation}
with the initial data
\begin{equation}\label{adv.11.ini}
 u(x,0)=\sin(x),\;\;\;x\in [0,2\pi].
\end{equation}
The rate of convergence of a numerical scheme, denoted by $r_l$, for the error $e_l$, is given by
 \begin{equation}
  r_l=\log(e_{l+1}/e_{l})/\log(N_{l+1}/N_l),\;\;l=0,1,2,\ldots,
 \end{equation}
where ${e}_l$ and $N_l$ denote the error and the number of mesh points, respectively. From Example \ref{ex.ordac}, we expect the rate of convergence, with respect to the space variable, of the CBSQI  
 and the QnBSQI schemes to be four and six, respectively. In order to achieve the rate numerically, we take 
 the time step $\Delta t=10^{-1}\Delta x^{1.5}$. The numerical solution is computed over the domain $[0,2\pi]$ at time $T = 1$. Table \ref{adv.cbsqi} and
   Table \ref{adv.qubsqi}, clearly shows that, as the number of mesh points increase, the numerical solution obtained with CBSQI and QnBSQI schemes converge to
    the exact solution with rate four and sixth, respectively.}
\end{example}
\begin{table}
\caption{Comparison of $L^{\infty}$-, $L^1$-, and $L^2$-errors
obtained with CBSQI scheme, in case of
\eqref{adv.11}-\eqref{adv.11.ini} with the exact solution over the
domain
$\Omega=[0,2\pi]$ at time $T=1$.}
\centering 
\small
\begin{tabular}{|l| c|c| c| c| c|  c| c| r|}
\hline 
 $N$ &$L^{\infty}$-error & order & $L^1$-error & order&$L^2$-error & order\\ 
\hline
20  &  3.188811e-04 & --     &0.001283 &--     &1.283344e-03 &--\\
\hline
40  &  2.020435e-05 & 3.980281     &8.097906e-05 &3.986215     &8.097906e-05 &3.986215\\
\hline
80  &  1.267336e-06 & 3.994795     &5.068834e-06 &3.997823     &5.068834e-06 &3.997823\\
\hline
160  &  7.925669e-08 & 3.999122     &3.170066e-07 &3.999069     &3.170066e-07 &3.999069\\
\hline
320  &  4.954097e-09 & 3.999839     &1.981728e-08 &3.999682     &1.981728e-08 &3.999682\\
\hline
\end{tabular}
\label{adv.cbsqi}
\end{table}
\begin{table}
\caption{Comparison of $L^{\infty}$-, $L^1$-, and $L^2$-errors
obtained with QnBSQI scheme, in case of
\eqref{adv.11}-\eqref{adv.11.ini} with the exact solution over the
domain
$\Omega=[0,2\pi]$ at time $T=1$.}
\centering 
\small
\begin{tabular}{|l| c|c| c| c| c|  c| c| r|}
\hline 
 $N$ &$L^{\infty}$-error & order & $L^1$-error & order&$L^2$-error & order\\ 
\hline
20  &  3.998450e-05 & --     &0.000026 &--     &2.840256e-05 &--\\
\hline
40  &  6.376781e-07 & 5.970469     &4.075987e-07 &5.974427     &4.522632e-07 &5.972714\\
\hline
80  &  1.013043e-08 & 5.976060     &6.453003e-09 &5.981035     &7.166100e-09 &5.979831\\
\hline
160  &  1.585622e-10 & 5.997504     &1.009394e-10 &5.998409     &1.121209e-10 &5.998061\\
\hline
320  &  2.633005e-12 & 5.912195     &1.673563e-12 &5.914423     &1.859182e-12 &5.914243\\
\hline
\end{tabular}
\label{adv.qubsqi}
\end{table}
We now show numerically that the CBSQI and QnBSQI schemes achieve the expected order of accuracy in the case of Burgers equation as long as the solution remains smooth.
\begin{example}\label{Burger.ord.ex}{\rm 

Consider the Burgers' equation 
\begin{equation}\label{bur.1}
 u_t+\Big(\frac{u^2}{2}\Big)_x=0.
\end{equation}
In order to show accuracy of the BSQI schemes, we consider
\eqref{bur.1} with the following initial data
\begin{equation}\label{bur.ini}
 u(x,0)=\sin(x),\;\;\;x\in [0, 2\pi]
\end{equation}
\begin{table}[t]
\caption{Comparison of $L^{\infty}$-, $L^1$-, and $L^2$-error
obtained with CBSQI scheme, in case of
Burgers' equation over the domain $\Omega=[0,2\pi]$ at time $T=0.5$.}
\centering 
\small
\begin{tabular}{|l| c|c| c| c| c|  c| c| r|}
\hline 
 $N$ &$L^{\infty}$-error & order & $L^1$-error & order&$L^2$-error & order\\ 
\hline
40  &  1.211681e-03 & --     &0.000161 &--     &3.518324e-04 &--\\
\hline
80  &  9.369829e-05 & 3.692843     &1.048627e-05 &3.937758     &2.465823e-05 &3.834747\\
\hline
160  &  6.419103e-06 & 3.867579     &6.624985e-07 &3.984440     &1.586015e-06 &3.958591\\
\hline
320  &  4.077475e-07 & 3.976624     &4.155074e-08 &3.994971     &9.982610e-08 &3.989845\\
\hline
640  &  2.553719e-08 & 3.997004     &2.598439e-09 &3.999157     &6.250058e-09 &3.997476\\
\hline
\end{tabular}
\label{table.CBSQI}
\end{table}
\begin{table}[t]
\caption{Comparison of $L^{\infty}$-, $L^1$-, and $L^2$-error obtained
with QnBSQI scheme, in case of
Burgers' equation over the domain $\Omega=[0,2\pi]$ at time $T=0.5$. }
\centering 
\small
\begin{tabular}{|l| c|c| c| c| c|  c| c| r|}
\hline 
 $N$ &$L^{\infty}$-error & order & $L^1$-error & order&$L^2$-error & order\\ 
\hline
40  &  3.802114e-04 & --     &0.000047 &--     &1.137541e-04 &--\\
\hline
80  &  1.388433e-05 & 4.775273     &9.976366e-07 &5.547309     &2.911152e-06 &5.288185\\
\hline
160  &  2.604317e-07 & 5.736409     &1.771954e-08 &5.815101     &5.314697e-08 &5.775459\\
\hline
320  &  4.307657e-09 & 5.917858     &2.849299e-10 &5.958591     &8.664328e-10 &5.938756\\
\hline
640  &  6.843121e-11 & 5.976105     &4.613797e-12 &5.948509     &1.370955e-11 &5.981835\\
\hline
\end{tabular}
\label{table.QnBSQI}
\end{table}
In Table \ref{table.CBSQI} and Table \ref{table.QnBSQI}, we compare
the $L^{\infty}$-, $L^1$-, and $L^2$-error and their rate of
convergence obtained using the CBSQI and the QnBSQI schemes, respectively. The numerical solution is
computed over the $[0,1]$ at time $T=0.5$, at which the solution remains smooth. 
In order to achieve the expected rate of
convergence, we choose the time step significantly small to make time
discretization error small. In both cases, we choose the time step $\Delta t=10^{-1}\Delta x^{1.5}$. From the tables, it is clear that the CBSQI scheme converges to the exact solution
with a rate more then four, whereas the QnBSQI scheme initial has order nearly four
which increased further to order more then six.
\qed}
\end{example}
\section{Hybrid Scheme}\label{hybrid.sec}

The CBSQI and the QnBSQI schemes developed in the previous section are stable with 3$^{\rm rd}$ order Runge-Kutta time discretization and also they achieve the expected order of accuracy when the solution remains smooth, as demonstrated in Example \ref{ex.ordacc.adv} and \ref{Burger.ord.ex}. But the main disadvantage of these schemes is that the numerical solution develops spurious oscillations in the case of non-smooth solution. In order to minimize the oscillations without much compromise in the order of accuracy, we devise an hybrid scheme, where we use the BSQI scheme in the region where the solution is smooth, whereas in the non-smooth region we propose to use WENO scheme. Thus, we propose the hybrid semi-discrete scheme for \eqref{scl.1d} to be
\begin{equation}
 \frac{du_j}{dt}=-\frac{1}{\Delta x}\Big( F_{j+\frac{1}{2}}^{\mbox{Hybrid}}-F_{j-\frac{1}{2}}^{\mbox{Hybrid}}\Big),
\end{equation}
where
\begin{equation}\label{hybrid.eq}
 F_{j+\frac{1}{2}}^{\mbox{\small Hybrid}}=\Phi_{j+\frac{1}{2}} F_{j+\frac{1}{2}}^{\mbox{WENO}}+(1-\Phi_{j+\frac{1}{2}}) F_{j+\frac{1}{2}}^{\mbox{BSQI}},
\end{equation}
with $\Phi$ being a smooth indicator, which is a constant over the $j^{\rm th}$ cell   {\it i.e}  $\Phi_{j-\frac{1}{2}}=\Phi_{j+\frac{1}{2}}=\Phi_j$. Using this in
 \eqref{hybrid.eq}, the above scheme can be written as
 \begin{align}\label{hybrid.semidisc}
   \frac{du_j}{dt}&=-\Phi_j\frac{1 }{\Delta x}\Big( F_{j+\frac{1}{2}}^{\mbox{WENO}}-F_{j-\frac{1}{2}}^{\mbox{WENO}}\Big)
   -(1-\Phi_j)\frac{1 }{\Delta x}\Big( F_{j+\frac{1}{2}}^{\mbox{BSQI}}-F_{j-\frac{1}{2}}^{\mbox{BSQI}}\Big).
 \end{align}
Here,  $F^{\mbox{BSQI}}$ denotes either CBSQI or QnBSQI flux given by \eqref{numflux.CBSQI} and \eqref{numflux.QnBSQI}, respectively. Whereas $F^{\mbox{WENO}}$ denotes the WENO flux of appropriate order.  More precisely, we conjugate WENO flux of order 3 with CBSQI and form a hybrid scheme with numerical flux
\begin{equation}\label{hybrid4.eq}
 F_{j+\frac{1}{2}}^{\mbox{\small Hybrid4}}=\Phi F_{j+\frac{1}{2}}^{\mbox{WENO3}}+(1-\Phi) F_{j+\frac{1}{2}}^{\mbox{CBSQI}},
\end{equation}
and we combine QnBSQI with WENO flux of order 5 to get 
 \begin{equation}\label{hybrid6.eq}
 F_{j+\frac{1}{2}}^{\mbox{\small Hybrid6}}=\Phi F_{j+\frac{1}{2}}^{\mbox{WENO5}}+(1-\Phi) F_{j+\frac{1}{2}}^{\mbox{QnBSQI}}.
\end{equation}
To be self-content, we quickly recall the WENO fluxes and for more details, we refer to Jiang and Shu \cite{jia-shu_96a}.

The flux $f$ can be decomposed into positive and negative parts as 
 \begin{equation}
  f(u)=f(u)^{+}+f(u)^{-},
 \end{equation}
 such that $df^{+}/du\geq 0$ and $df^{-}/du\leq 0$. For instance, these can be taken as 
 $$f(u)^{\pm} = f(u) \pm \alpha u,$$
 where $\alpha=\max|f'(u)|$.
We recall that the stencil of WENO3 is $S_0=\{x_{j-1},x_j\}$ and $S_1=\{x_j,x_{j+1}\}$, and the corresponding numerical flux is written as
$$ F_{j+\frac{1}{2}}^{\mbox{WENO3}} =  (F_{j+\frac{1}{2}}^{\mbox{WENO3}})^+ + ( F_{j+\frac{1}{2}}^{\mbox{WENO3}})^-$$
where the positive part is given by
\begin{eqnarray*}
\left(F_{j+\frac{1}{2}}^{\mbox{WENO3}}\right)^+=\omega_0 q_0^1+\omega_1 q_1^1,
\end{eqnarray*}
with
\begin{eqnarray*}
q_0^1=-\frac{1}{2}f_{j-1}^{+}+\frac{3}{2}f_j^{+},\;\;\; q_1^1=\frac{1}{2}f^+_{j}+\frac{1}{2}f_{j+1}^{+},~~\omega_k=\frac{\alpha_k}{\displaystyle{\sum_{l=0}^1 \alpha_l}},\;\;\;k=0,1,\\
\alpha_k=\frac{C^1_k}{(\epsilon+IS_k)^2},\;\;k=0,1,~~ IS_0=(f_{j}^{+}-f_{j-1}^{+})^2,\;\;\;IS_1=(f_{j+1}^{+}-f_{j}^{+})^2,
\end{eqnarray*}
and $C^1_0=1/3$ and $C^1_1=2/3$. The negative part can be written in a similar way. 

The stencil of WENO5 is 
\begin{align*}
  S_0=\{x_{j-2},x_{j-1},x_{j}\}, S_1=\{x_{j-1},x_{j},x_{j+1}\},S_2=\{x_{j},x_{j+1},x_{j+2}\},
 \end{align*}
 and the positive part of the numerical flux is written as
\begin{equation*}
\left(F_{j+\frac{1}{2}}^{\mbox{WENO5}}\right)^+=\displaystyle{\sum_{k=0}^2 \omega_kq_k^2},
 \end{equation*}
 where
 \begin{align*}
  q_0^2=\frac{1}{3}f_{j-2}^{+}-\frac{7}{6}f_{j-1}^{+}+\frac{11}{6}f_j^{+},~
  q_1^2=-\frac{1}{6}f_{j-1}^{+}+\frac{5}{6}f_{j}^{+}+\frac{1}{3}f_{j+1}^{+},~
  q_2^2=\frac{1}{3}f_{j}^{+}+\frac{5}{6}f_{j+1}^{+}-\frac{1}{6}f_{j+2}^{+},
 \end{align*}
and
\begin{equation}
 \omega_k=\frac{\alpha_k}{\displaystyle{\sum_{l=0}^2 \alpha_l}},~~ \alpha_k=\frac{C^2_k}{(\epsilon+IS_k)^2}, \;\;\;k=0,1,2.
\end{equation} 
 Here the smooth indicators are given by
 \begin{align*}
  IS_0&=\frac{13}{12}(f^{+}_{j-2}-2f^{+}_{j-1}+f^{+}_{j})^2+\frac{1}{4}(f^{+}_{j-2}-4f^{+}_{j-1}+3f^{+}_{j})^2,\\
  IS_1&=\frac{13}{12}(f^{+}_{j-1}-2f^{+}_{j}+f^{+}_{j+1})^2+\frac{1}{4}(f^{+}_{j-1}-f^{+}_{j+1})^2,\\
  IS_2&=\frac{13}{12}(f^{+}_{j}-2f^{+}_{j+1}+f^{+}_{j+2})^2+\frac{1}{4}(3f^{+}_{j}-4f^{+}_{j+1}+f^{+}_{j+2})^2,
 \end{align*}
and $C_1 = 1/10, C_2 = 6/10$, and $C_3 = 3/10$.

Thus, the hybrid scheme switches between WENO and BSQI scheme, depending upon the smooth indicator. Once the shock region is detected, the 
  hybrid scheme switches to the expensive shock capturing WENO scheme, otherwise the solution is computed using the BSQI scheme.
  
For the smooth indicator, we use the weak local truncation error (WLTE) introduced by Karni and Kurganov \cite{kar-kur_05a}. The WLTE is derived on 
the fact that a weak solution of \eqref{scl.1d}
satisfies the integral form given by (see Godlewski and Raviart \cite{god-rav_91a} )
 \begin{equation}\label{wlr}
  E(u,\phi):=-\int_0^T\int_\mathbb{R}\{ u(x,t)\phi_t(x,t)+f(u)\phi_{x}(x,t)\}dxdt+\int_\mathbb{R}u(x,0)\phi(x,0) dx=0,
 \end{equation}
 for all $\phi(x,t)\in C_0^1(\mathbb{R}\times (0,T])$. The variation in the value of $E$ from smooth to discontinuous regions and vice versa for computed solution $u^{\Delta}$ can be taken as a measure of smooth indicator.  The $E(u^{\Delta},\phi)$ is referred as {\it weak truncation error} for $u^{\Delta}$ with respect to $\phi$. In practice, the calculation of 
weak truncation error is seem to be a difficult task, since $\phi$ is a general test function. To overcome this difficulty Kurganov and his co-workers
(see \cite{kur-liu_12a,kar-kur_05a}) used test function with B-splines as 
 \begin{equation}\label{Bspline.testf}
  \phi_{j}^{n}(x,t)=B_{j}(x)B^{n}(t),
 \end{equation}
 where 
 \begin{equation}
  B_{j}(x) = 
 \left\{\begin{array}{ll}
\dfrac{1}{2}\left(\dfrac{x-x_{j-3/2}}{\Delta x}\right)^2 &\mbox{ if } x_{j-3/2}<x\le x_{j-1/2},\medskip\\
\dfrac{3}{4}-\left(\dfrac{x-x_{j}}{\Delta x}\right)^2&\mbox{ if } x_{j-1/2}<x\le x_{j+1/2},\medskip\\
 \dfrac{1}{2}\left(\dfrac{x-x_{j+3/2}}{\Delta x}\right)^2  &\mbox{ if } x_{j+1/2}<x\le x_{j+3/2},\medskip\\
0&\mbox{otherwise},
\end{array}\right.
\end{equation}
and

\begin{equation}
  B^{n}(t) = 
 \left\{\begin{array}{ll}
\left(\dfrac{t-t^{n-3/2}}{\Delta t}\right) &\mbox{ if } t^{n-3/2}<t\le t^{n-1/2},\medskip\\
 \left(\dfrac{t^{n+1/2}-t}{\Delta t}\right)&\mbox{ if } t^{n-1/2}<t\le t^{n+1/2},\medskip\\
0 &\mbox{otherwise}
\end{array}\right.
\end{equation}
are the quadratic and the linear B-splines with the localized supports of size $|\mbox{supp}(B_{j})|=3 \Delta x$ and  $|\mbox{supp}(B^{n})|=2 \Delta t$, respectively.  Putting \eqref{Bspline.testf} in \eqref{wlr}, we can arrive at 
\begin{multline}\label{wlte}
 E^{n}_{j}=\frac{1}{6}\Big[ u^n_{j+1}-u^{n-1}_{j+1}+4(u^n_{j}-u^{n-1}_{j})
 +u^n_{j-1}-u^{n-1}_{j-1}\Big]\Delta x \\
 + \frac{1}{4}\Big[ f(u^n_{j+1})-f(u^n_{j-1}) +f(u^{n-1}_{j+1})-f(u^{n-1}_{j-1})\Big ]\Delta t.
 \end{multline}
Karni and Kurganov \cite{kar-kur_05a} established an estimate for the WLTE \eqref{wlte} and
remarked that the WLTE bound can be converted into $L^{\infty}_{{loc}}$ error bound in the smooth region. Although their result is valid for one dimension scalar
  problem, they showed numerical that the WLTE behavior is similar in the case of 1D system of equations.  With the aid of these estimates and further numerical experiments, Kurganov and Liu \cite{kur-liu_12a} proposed the following estimate and used it  as a smooth indicator:
 \begin{equation*}\label{wlre.est}
  \Arrowvert E^n_j\Arrowvert_{\infty}\approx \left\{\begin{array}{ll}
\Delta , \mbox{near the shock,}\medskip\\
\Delta^{\alpha} , \mbox{near the contact wave},\;\;1<\alpha \leq 2,\medskip\\
 \Delta^\beta, \mbox{in the smooth region,} \medskip\\
\end{array}\right.
 \end{equation*}
where $\Delta =\max(\Delta x,\Delta t)$, $\beta=\min\{r+2,4\}$ and $r$ is the accuracy of the numerical scheme.
In our hybrid scheme \eqref{hybrid.semidisc}, we propose to use the smooth indicator defined by
\begin{equation}\label{sind.est}
\Phi(x_j)=\left\{\begin{array}{ll}
                1 ~~~~  | E_j|>K\Delta x^4, \\
                0 ~~~~ ~~ \text{otherwise},
                \end{array}\right.
\end{equation}
where $K$ is some empirically chosen nonnegative real number. In general, we take $K=1/\Delta x$ in our numerical experiments. To make smooth transition between smooth and rough 
 regions, we extend the rough region to the neighboring node points of $x_j$. That is, once the discontinuity is detected in a vicinity of $x_j$, which is the case if $\Phi(x_j)=1$, then we take
 \begin{equation}
  \Phi(x_{j\pm a})=1, \;\;a=1,2,\ldots,M,
 \end{equation}
 where $M$ is again an empirically chosen positive integer.
 In our numerical computation, we take $M=2$. The smooth transition can also be made using a smaller value of $K$, but that lead to 
  over estimate of discontinuous regions.
\section{Numerical Experiments}\label{numericalhybrid}

This section comprises the numerical implementation of hybrid BSQI-WENO schemes for linear and non-linear problems. The BSQI schemes are oscillatory in nature at shock positions,  whereas the hybrid BSQI-WENO schemes developed in the previous section are non-oscillatory because of the usage of the WENO scheme in the non-smooth regions and the high order accuracy is maintained in the smooth regions using BSQI approximation. The developed algorithm works for any degree of BSQI, but to demonstrate the idea we use CBSQI and QnBSQI for space approximation in conjugation with WENO3 and WENO5, respectively. We study the accuracy of the CBSQI-WENO3 and the QnBSQI-WENO5 schemes in terms of rate of convergence. Along with the accuracy, we also demonstrate the efficiency of the proposed hybrid BSQI-WENO. 

We start the numerical experiment with smooth initial data.
\begin{figure}
 \includegraphics[width=8.3cm]{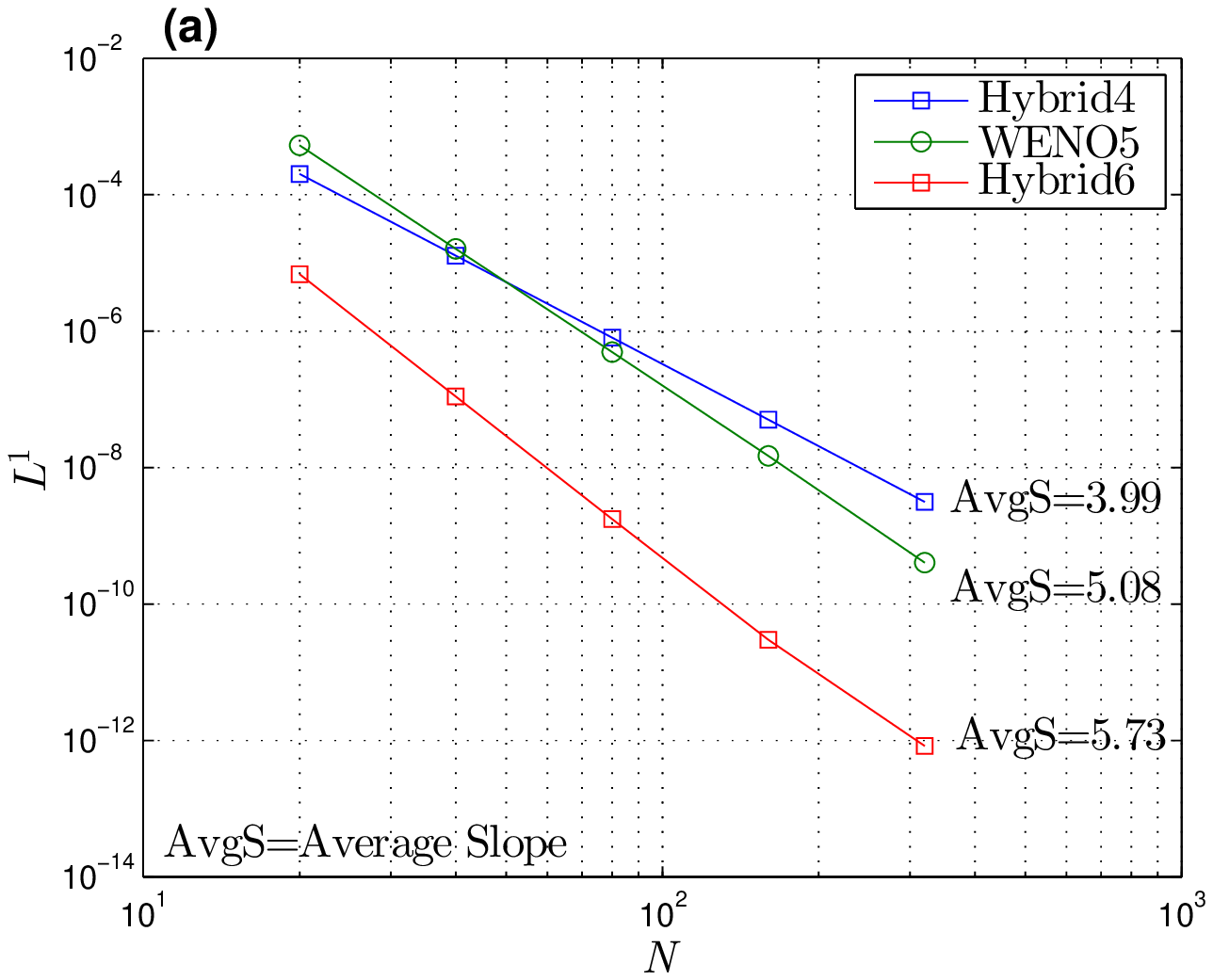}
 \includegraphics[width=8.3cm]{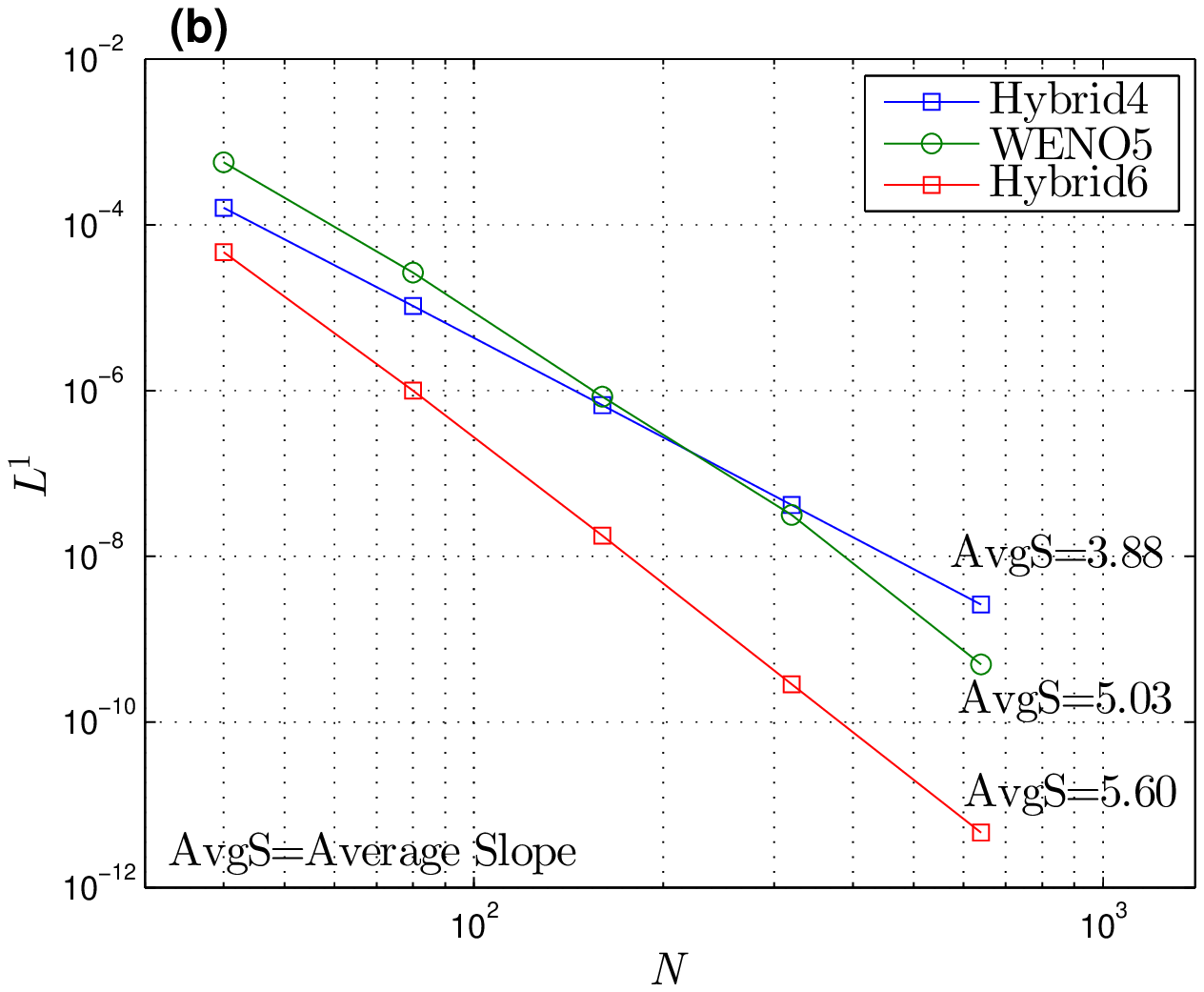}
 \caption{Comparison of the Hybrid4 and the Hybrid6 schemes with the WENO5 scheme in terms of the $L^1$-error (a) for the problem
\eqref{adv.11}-\eqref{adv.11.ini} at time $T=1$, and
 (b) for the problem \eqref{bur.1}-\eqref{bur.ini} at time $T=0.5$.}
 \label{hybrid_con}
\end{figure}
\begin{example}[Smooth Initial Data]{\rm

In order to measure the accuracy and convergence rate of the hybrid algorithm, we consider the Example 3.1 and Example 3.2. In  both cases, we take $K=1$ in smooth indicator. We fixed the time step $\Delta t=10^{-1}\Delta x^{1.5}$ and increase the mesh points  in space. In Figure \ref{hybrid_con}(a), we compare the accuracy of the Hybrid4 and the Hybrid6 schemes with the WENO5 scheme in case of linear advection equation at time $T=1$. Both the Hybrid4 and the Hybrid6 schemes achieve their expected convergence rate. In terms of accuracy, we observe that the Hybrid6 scheme is better than both the WENO5 and the Hybrid4 schemes. Similar behavior is observed from Figure \ref{hybrid_con}(b) in the case of Burgers' equation at time $T=0.5$. }
\end{example}

As a next step, we perform the numerical experiment in the case of the linear advection equation with discontinuous initial data.
\begin{figure}
 \includegraphics[width=0.52\textwidth]{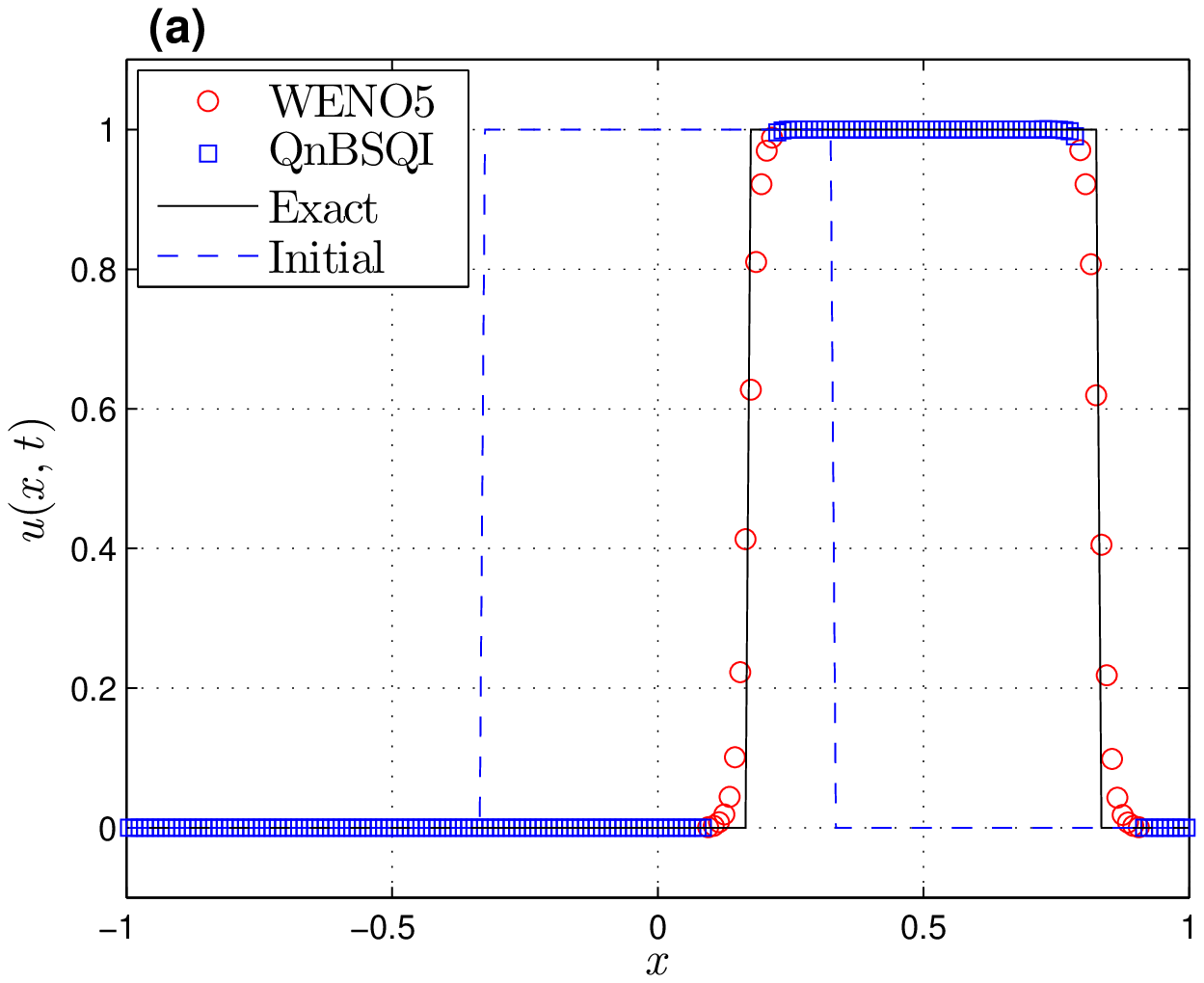}
 \includegraphics[width=0.52\textwidth]{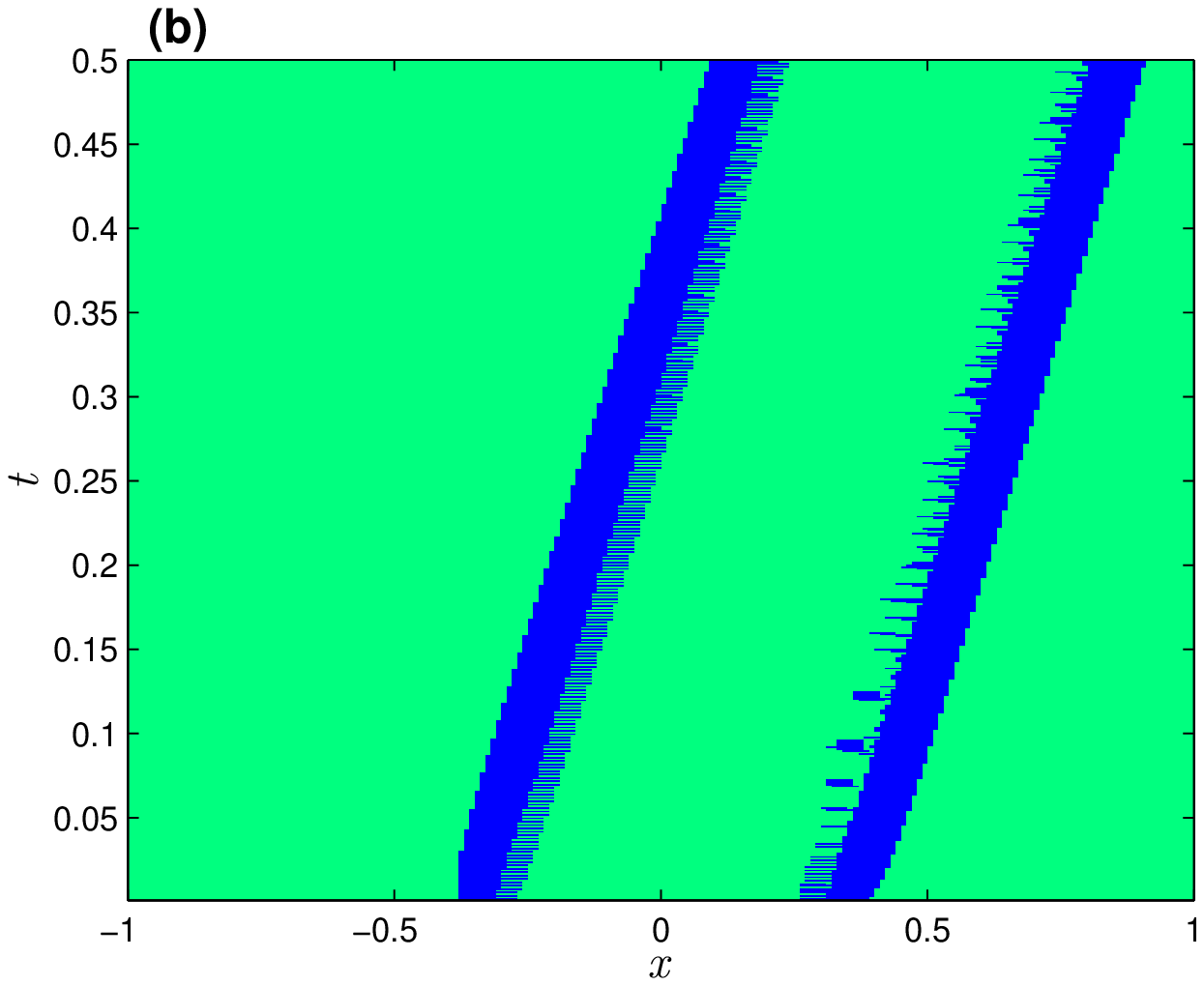}
 \caption{(a) Comparison of the numerical solution obtained using the Hybrid6 scheme \eqref{hybrid6.eq}
 with the exact solution at time $t=.5$, (b) Plot of the smooth indicator \eqref{sind.est} for the time ranging in $[0, 0.5]$. Here blue color indicates the regions where WENO5 scheme is used ({\it i.e.} $\Phi=1$) and green color indicates the region where QnBSQI scheme is used. }
 \label{adv.fig}
\end{figure}
\begin{example}[Linear advection equation]{\rm

Consider the linear advection equation \eqref{adv.11} with initial data
  \begin{equation}\label{buk.ini1}
u(x,0)=\left\{\begin{array}{ll}
                1 ~~~~   \text{if} ~~ |x|\leq 1/3, \\
                0 ~~~~  \text{otherwise},
                \end{array}\right.
\end{equation}
which has jump discontinuities at points $x=-1/3,1/3$ and moves towards right as time progress. The CBSQI and QnBSQI schemes capture the correct position of the discontinuities but with oscillations in the profile (not shown here).  We now compute the solution using the Hybrid6 scheme whose numerical flux is given by \eqref{hybrid6.eq}.
 In Figure \ref{adv.fig}(a), we compare the numerical solution with the exact solution at time $t=0.5$ along with the initial data. The number of mesh points is taken to be $200$ and the
 CFL number is $0.4$. In the figure, the symbol `$\degree$' represents the region where the WENO5 is used whereas the symbol `$\square$' represents the region where QnCBSQI is used.  We observe that the discontinuities are captured without oscillations. In Figure \ref{adv.fig}(b), the plot of smooth indicator over the $(x,t)$ plane is shown. It is clear from the figure that the smooth indicator moves along the characteristic direction and detect the discontinuity positions correctly.  Also, we observe that the solution is computed with WENO5 only at discontinuity positions, whereas in the smooth regions QnBSQI  is used to compute the solution.\qed
}\end{example}

We now test the Hybrid6 scheme in the case of some nonlinear scalar conservation laws.  We first start with the inviscid Burgers equation, which is a common test case for convex flux.
 \begin{figure}[]
 \includegraphics[width=0.5\textwidth]{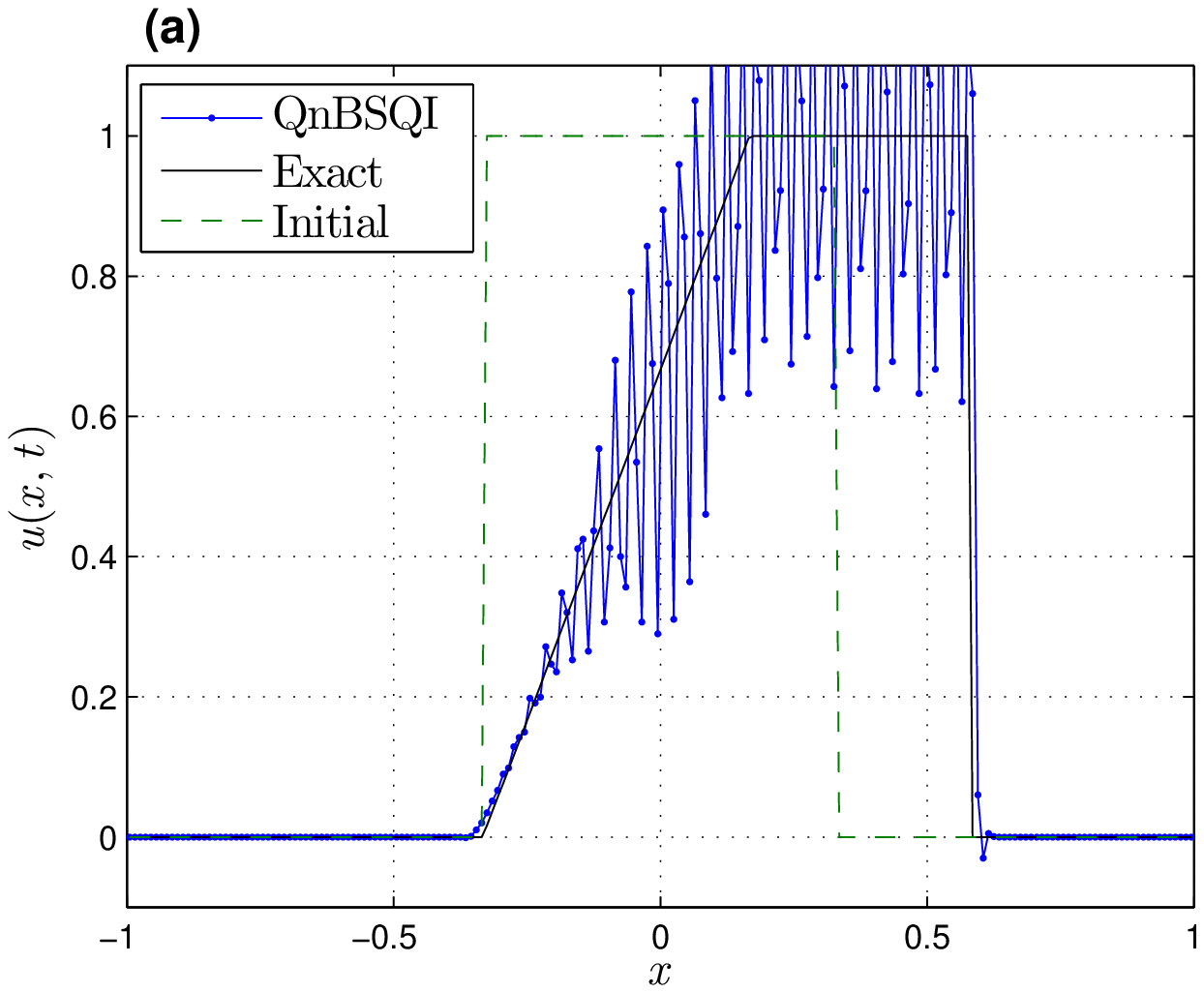}
 \includegraphics[width=0.5\textwidth]{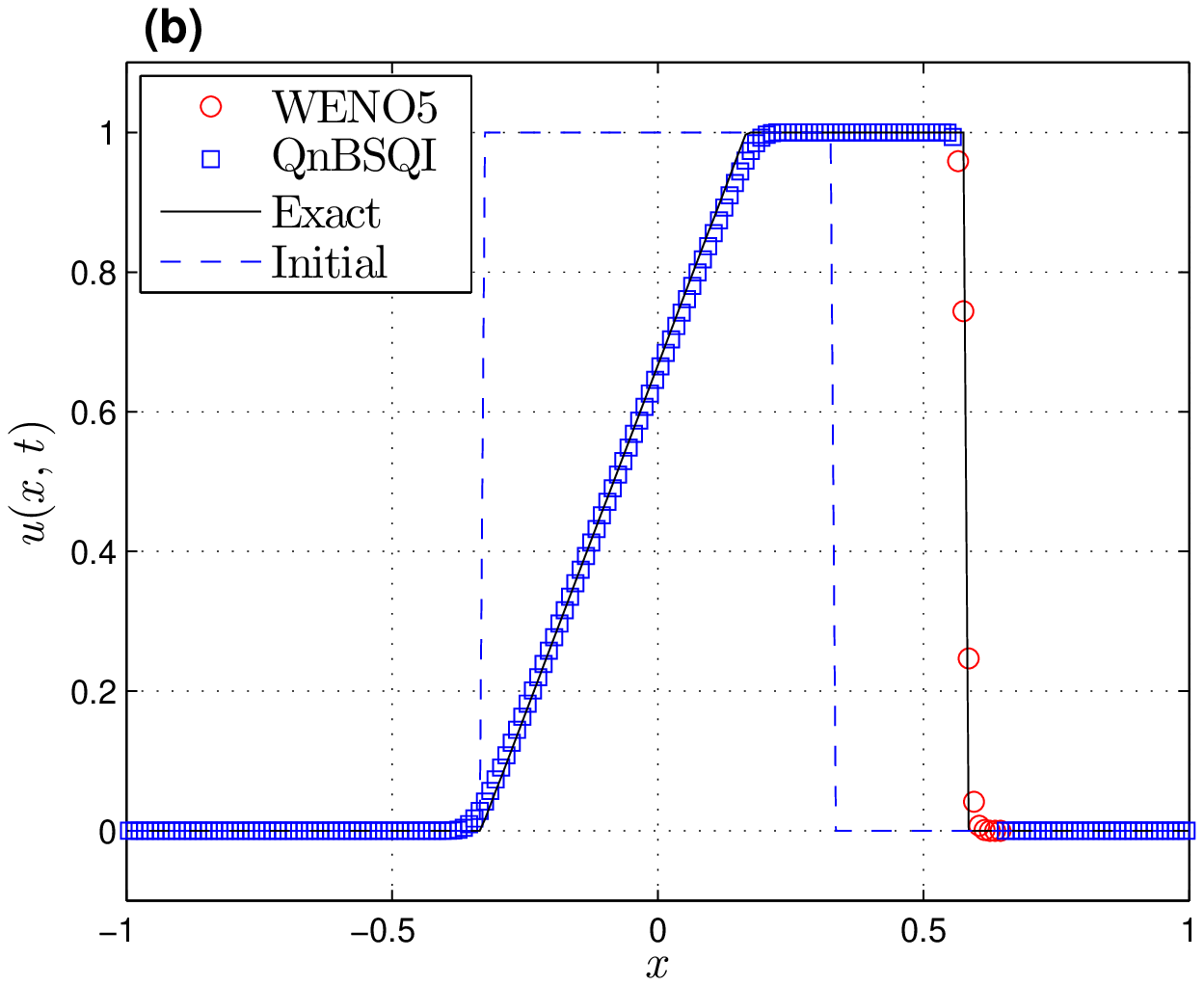}
 \includegraphics[width=0.5\textwidth]{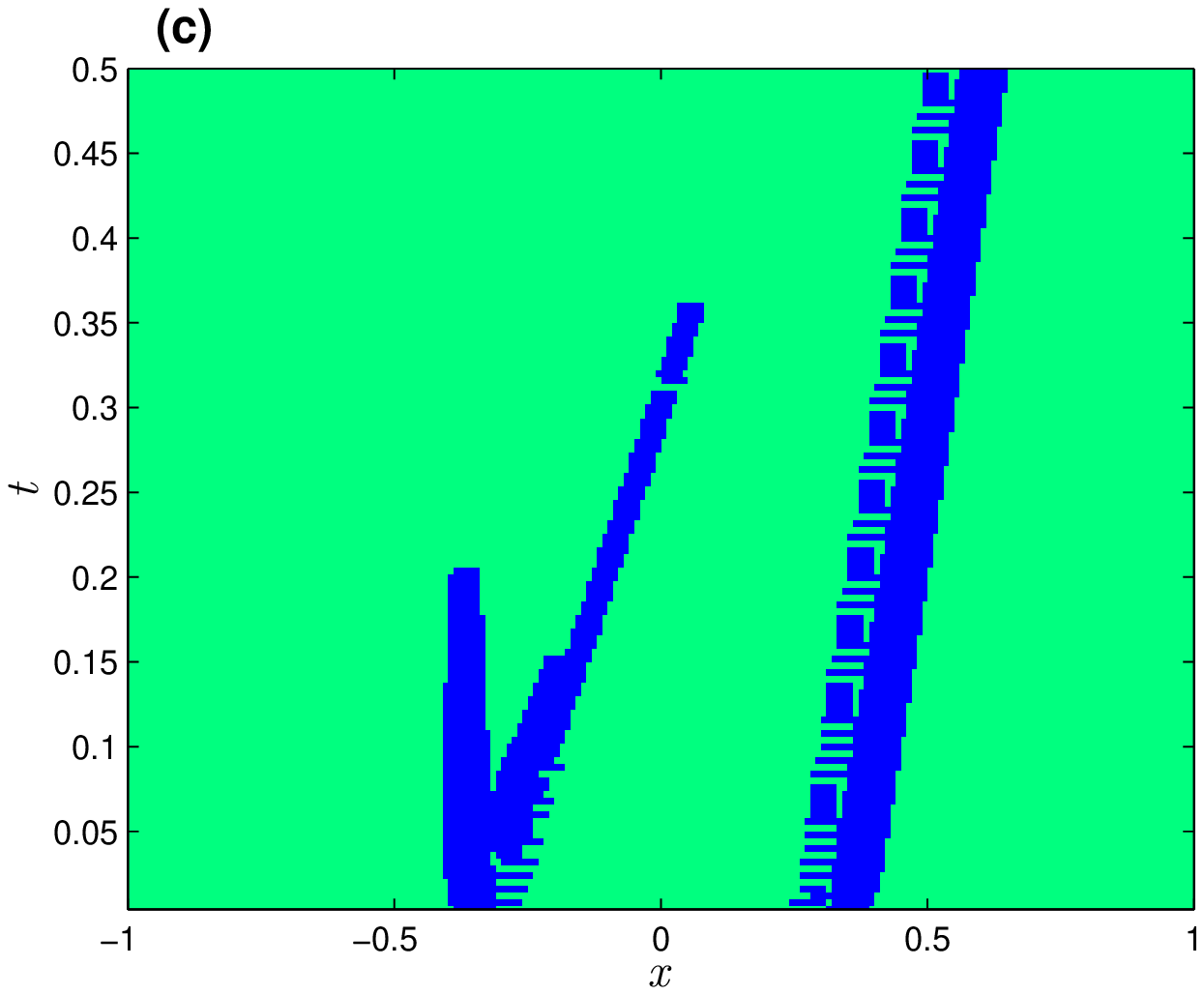}
 \includegraphics[width=0.5\textwidth]{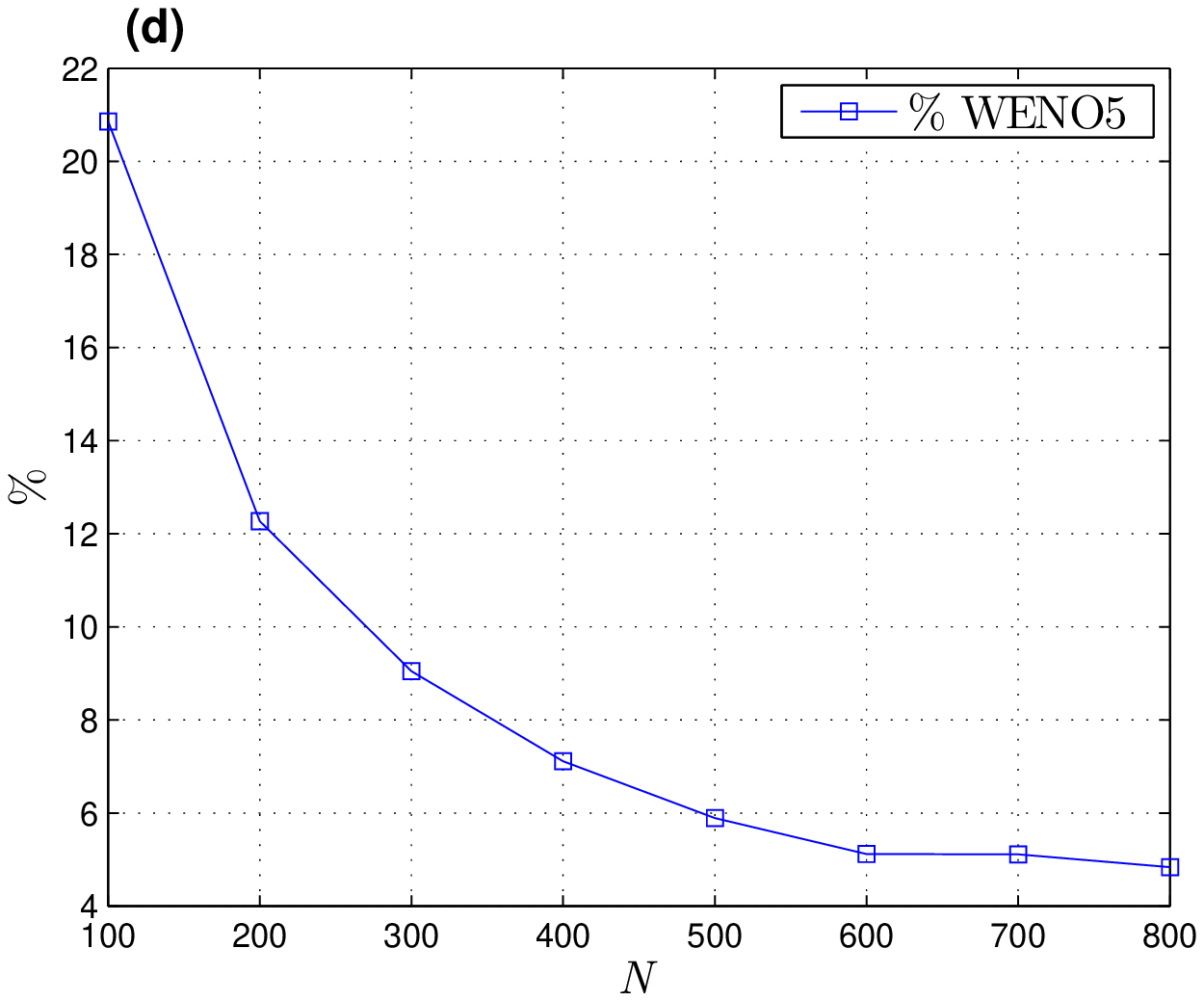}
 \caption{(a) Comparison of the numerical solution for the Burgers' equation with initial condition \eqref{buk.ini1} obtained by QnBSQI with the reference solution, 
at time $t=0.5$, (b) Comparison of the numerical solution of the same problem obtained using Hybrid6 given by \eqref{hybrid6.eq},
(c) Plot of smooth indicator for time interval $[0, 0.5]$, where the blue region indicates the usage of WENO5 and the green region indicates the usage of QnBSQI,
(d) Percentage of WENO5 scheme used in Hybrid6  to compute solution with increase in number of mesh points.}
 \label{bur.1d}
\end{figure}

\begin{example}[Nonlinear equation with convex flux - Burgers equation]{\rm

To show the shock capturing ability of the Hybrid6 scheme \eqref{hybrid6.eq}, we consider the Burgers' equation \eqref{bur.1} with discontinuous initial data \eqref{buk.ini1}. The numerical solution is computed at time $t=0.5$ over the domain $[-1,1]$ with $200$ mesh points. The CFL number is taken to be $0.4$. The solution contains a rarefaction wave and a shock moving towards right. As expected the QnBSQI scheme is oscillatory at the shock position, which spreads over the smooth region of the solution as the time progresses, which can be observed from Figure \ref{bur.1d}(a). 

Figure \ref{bur.1d}(b) depicts that the numerical solution at time $t=0.5$ obtained using the Hybrid6 scheme \eqref{hybrid6.eq}.  Here, we observe that the scheme captures shock and rarefaction accurately without any oscillation. The smooth indicator over the $(x,t)$-plane is shown in Figure \ref{bur.1d}(c), where we observe that the indicator moves along the discontinuity. The left state of the initial discontinuity at the point $x=-1/3$  is connected with the right state by the rarefaction fan, where the singularity at the lower part of the rarefaction wave is smoothen after time $t=0.25$ due to numerical dissipation, which is well illustrated by the smooth indicator.

In Figure \ref{bur.1d}(d), we show the plot of the percentage of the computational domain in which WENO5 is used verses the number of grid points.  We kept the time fixed at $t=0.5$ and increased the number of mesh points and observed that the percentage of the region where WENO5 scheme is used has been reduced drastically from $21\%$ to $5\%$ when the number of mesh points is increased from $100$ to $800$. The low percentage of usage of the WENO scheme for large number of mesh points, makes the hybrid BSQI-WENO schemes efficient.\qed
}
\end{example}
We now turn our attention to \eqref{scl.1d} with non-convex fluxs.  Since $f''(u)$ vanishes at one or more points, \eqref{scl.1d} fails to be genuinely nonlinear which makes the solution more complicated as it may involve the formation of the composite waves.  The composite waves are generally the combination of joined rarefaction and shock waves, for more detail see LeVeque \cite{lev_02a}. Many popular schemes fail to approximate the entropy solution of \eqref{scl.1d} with non-convex flux, for instance semi-discrete central-upwind schemes (see Kurganov {\it et al.} \cite{kur-etal_07a}). Here we test the Hybrid6 scheme  \eqref{hybrid6.eq} in case of the Buckley-Leverett problem.
 \begin{figure}
  \includegraphics[width=0.5\textwidth]{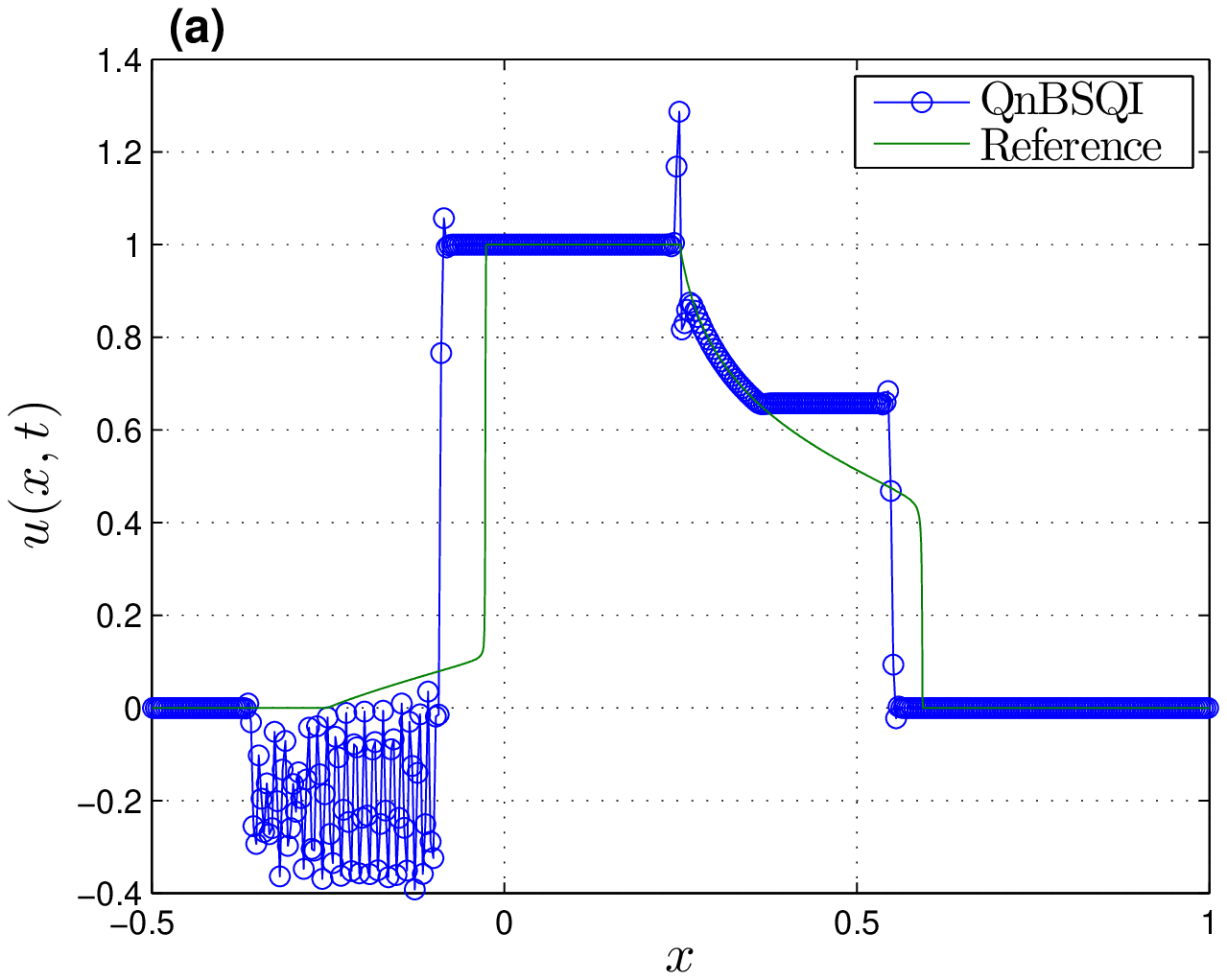}
 \includegraphics[width=0.5\textwidth]{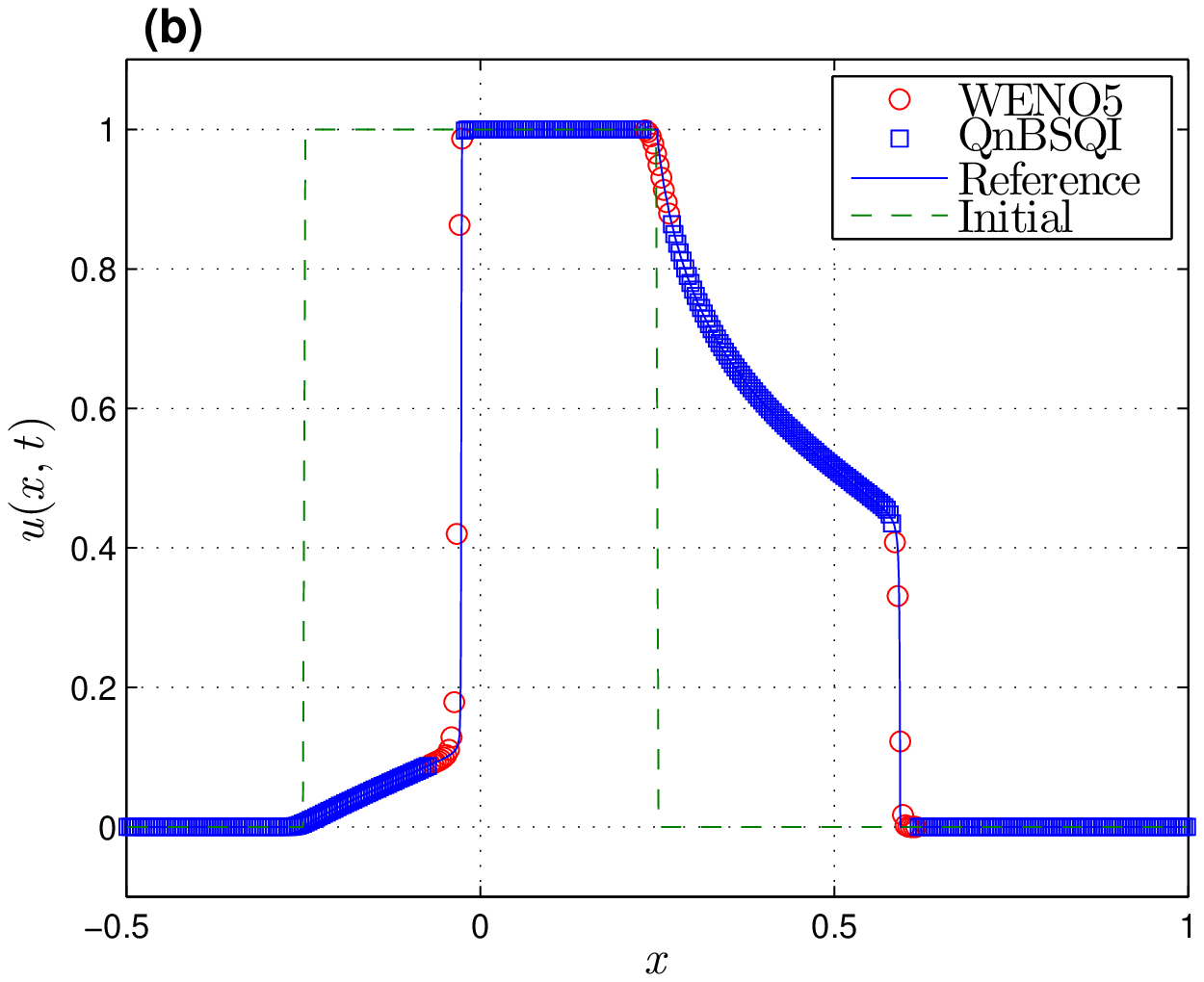}
 \includegraphics[width=0.5\textwidth]{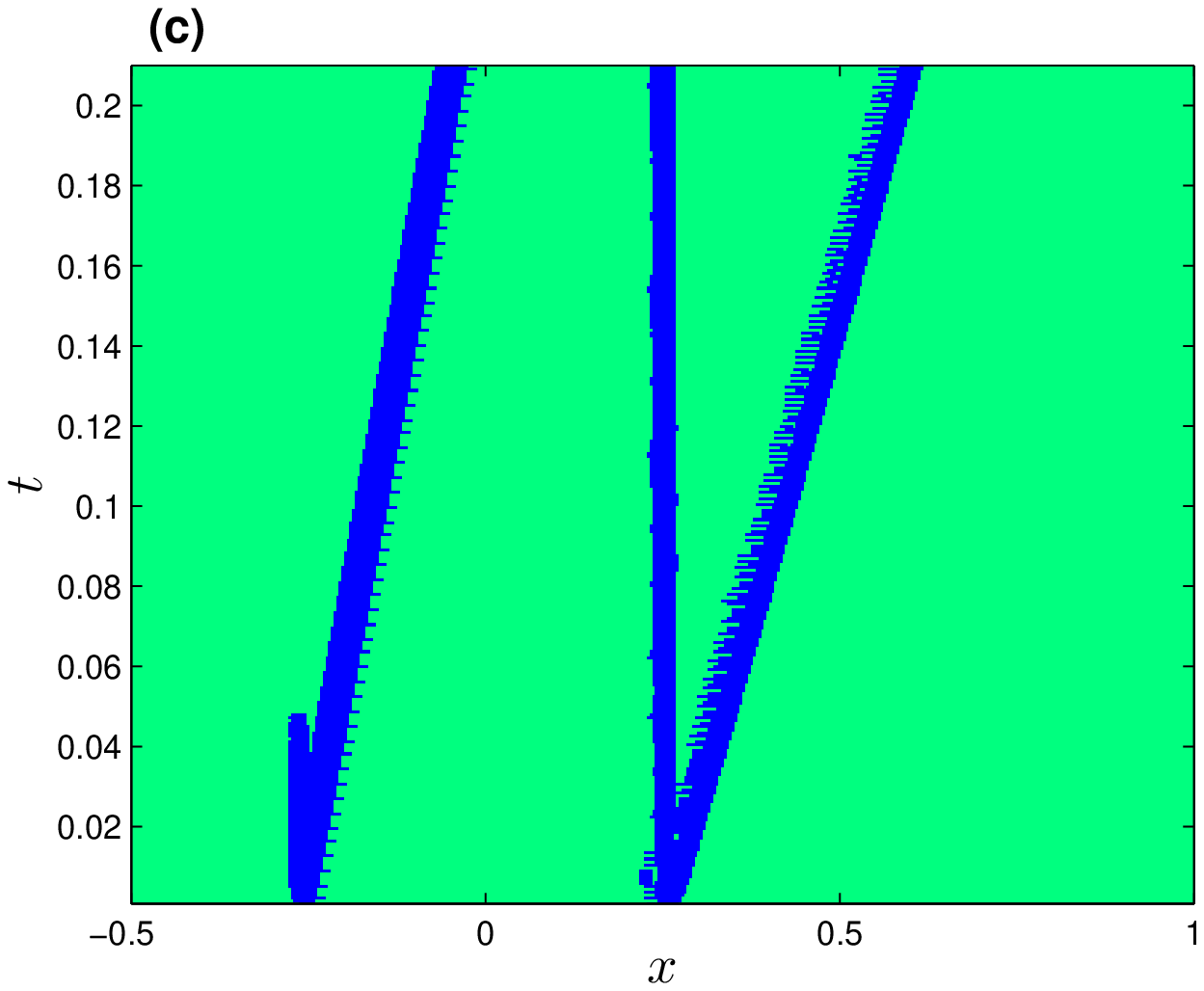}
 \includegraphics[width=0.5\textwidth]{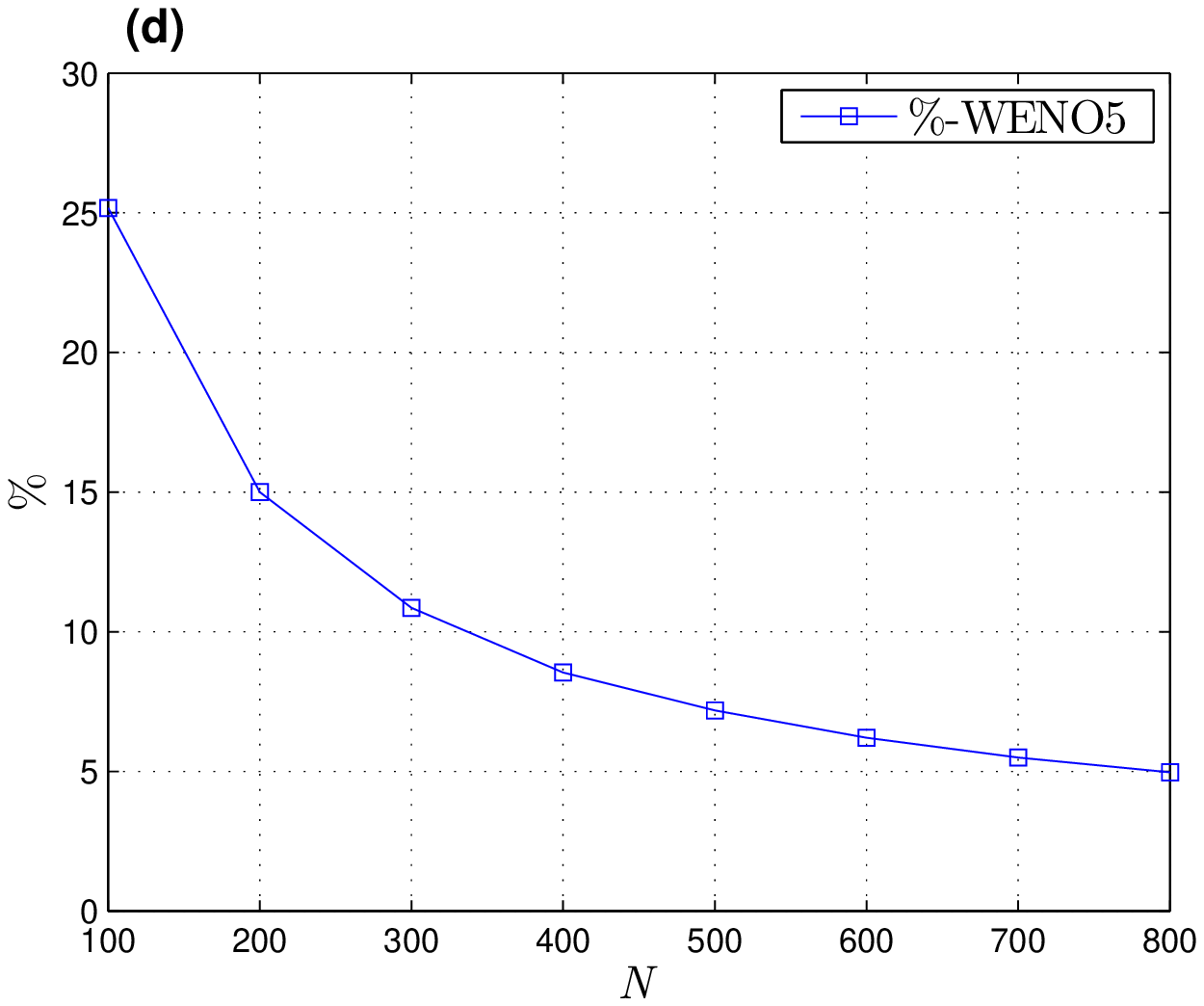}
 \caption{a) Comparison of numerical solution obtained by QnBSQI with reference solution, 
at time $t=0.21$ with initial data, (b) Hybrid6   (c) Plot of smooth indicator for time interval $[0, 0.21]$,
(d) Percentage of WENO5 scheme is used in Hybrid6  to compute solution with increase in number of mesh points.}
 \label{buc.1d}
\end{figure}

\begin{example}[Non-convex flux - Buckley-Leverett equation]{\rm

Consider the conservation law \eqref{scl.1d}  with non-convex flux, called the {\it Buckley-Leverett equation}, 
 \begin{equation}\label{BL.flux}
  f(u)=\frac{u^2}{(u^2+(1-u^2))},
 \end{equation}
along with the initial data \eqref{buk.ini1}. In Figure \ref{buc.1d}(a), the numerical solution obtained using QnBSQI scheme at time $t=0.21$ is compared with the reference solution over the domain $[-1,1]$ with $\Delta x=1/400$ and CFL is taken to be $0.2$. In this figure, we observe that the QnBSQI scheme fails to produce non-oscillatory result. The rarefaction fan is missing and the numerical solution is oscillatory at the rarefaction position, which makes the solution to moves with different speed. As a result, the QnBSQI scheme fails to capture the correct positions of the shocks. In comparison to QnBSQI, the Hybrid6 scheme \eqref{hybrid6.eq} produces non-oscillatory solution with correct positions of shocks as is apparent from Figure \ref{buc.1d}(b) and thus converges to the entropy solution. In Figure \ref{buc.1d}(c), we show the smooth indicator on $(x,t)$-plane and the percentage of the usage of WENO5 scheme is depicted in Figure \ref{buc.1d}(d).  As observed in the case of Burgers equation, here also we see that the usage of WENO5 is approximately $5\%$ when the number of mesh points approaches $800$. This shows that, the Hybrid6 scheme captures the entropy solution with less computational complexity.
\qed}
\end{example}
We continue with the scalar conservation law \eqref{scl.1d} with another non-convex flux.
\begin{example}[Another non-convex flux]{\rm

Another test case for non-convex flux is (see for instance Kurganov {\it et al.} \cite{kur-etal_07a})
\begin{equation}\label{buk.ini11}
f(u)=\left\{\begin{array}{ll}
                \dfrac{u(1-u)}{4} ~~~~& u<\dfrac{1}{2}, \medskip\\
               \dfrac{1}{2}u^2-\dfrac{1}{4}u+\dfrac{3}{16} ~~~~& u \geq \dfrac{1}{2},
                \end{array}\right.
\end{equation}
which belongs to $C^1$, whose convexity changes at $u=1/2$. We consider the following initial data 
 \begin{equation}\label{buk.ini12}
u(x,0)=\left\{\begin{array}{ll}
                1 ~~~~    x\leq 0.25, \\
                0 ~~~~    x>0.25.
                \end{array}\right.
\end{equation}
The exact solution is given by
 \begin{equation}\label{buk.ini43}
u(x,t)=\left\{\begin{array}{ll}
                0 & ~~~x<\dfrac{((\sqrt{6}-2)t+1)}{4}, \medskip\\
                \dfrac{(x-0.25)}{t}+\dfrac{1}{2} &~~~\dfrac{((\sqrt{6}-2)t+1)}{4}<x<\dfrac{2t+1}{4},\medskip\\
                1&~~~x>\dfrac{2t+1}{4}.
                \end{array}\right.
\end{equation}
Capturing the entropy solution for the above initial-value problem using higher order methods is a challenging problem. Kurganov {\it et al.} \cite{kur-etal_07a} illustrated that the numerical solutions of many schemes does not converge to the unique entropy solution or convergence is too slow that, it requires very fine meshes to achieve a required accuracy. They developed an adaptive dissipative minmod reconstruction in the domain where the flux changes the sign and for the other parts of domain, the fifth order WENO reconstruction is used.  In our proposed Hybrid6 scheme, without making any modifications, the algorithm follow the same procedure as in the convex case. In Figure \ref{buc.ex1}(a), we compare the solution obtained using Hybrid6 scheme with exact solution at time $t=1$, where we have taken the CFL number equals 0.2 and the number of mesh points is 200.   From this comparison result, we observe that the numerical solution seems to be converging to the exact entropy solution \eqref{buk.ini43}.  From the Figure \ref{buc.ex1}, it is clear that, WENO5 is used only at shock position and rarefaction tip whereas the QnBSQI approximates the solution in smooth region at $t=1$. The scheme is efficient as very low percentage of WENO scheme is used. 
\begin{figure}
 \includegraphics[width=0.5\textwidth]{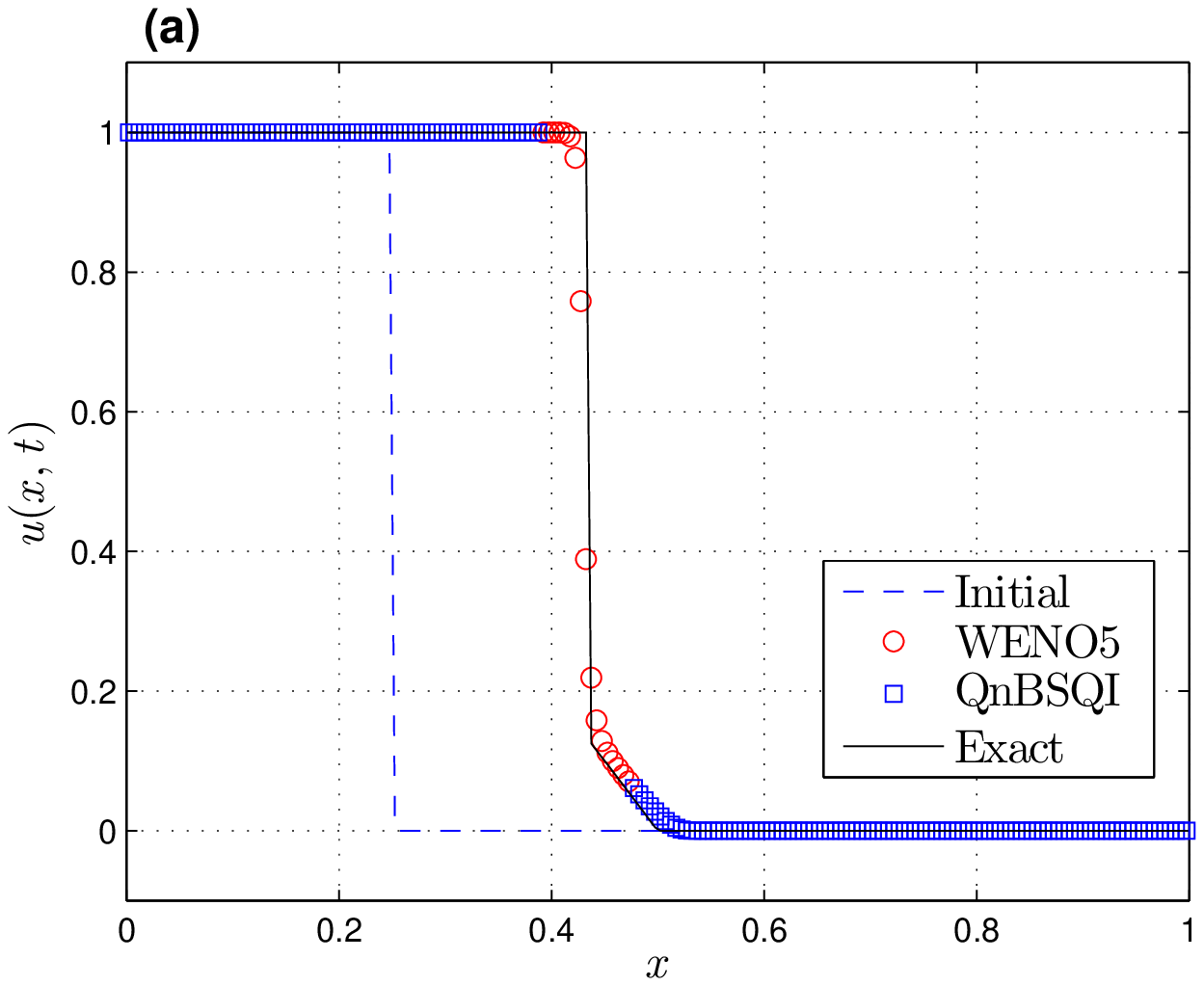}
 \includegraphics[width=0.5\textwidth]{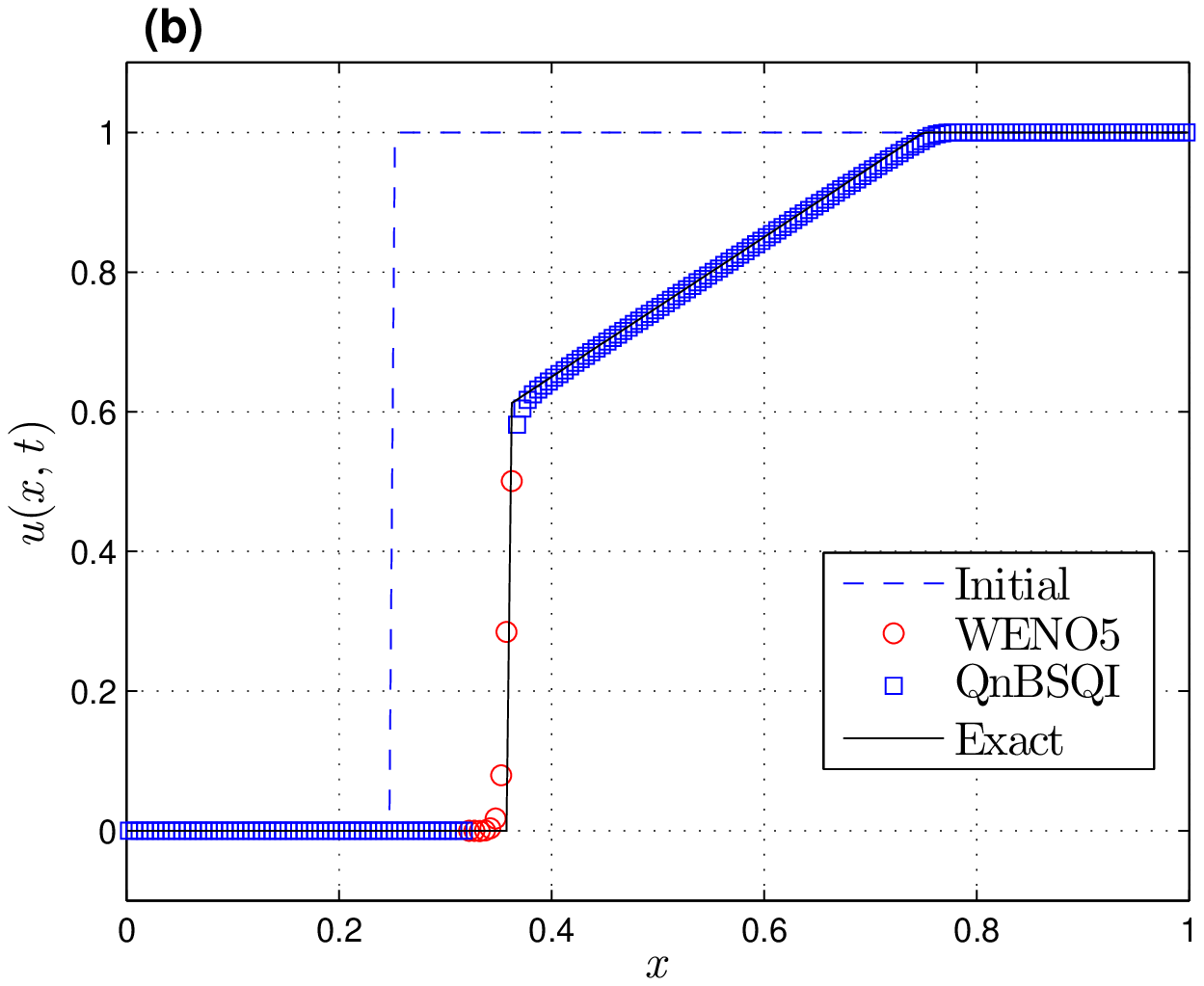}
 \caption{(a) Comparison of the numerical solution for the problem \eqref{buk.ini11}-\eqref{buk.ini12} obtained by Hybrid6 with exact solution \eqref{buk.ini43} at time $t=1$,
 (b)  Comparison of the numerical solution for the problem \eqref{buk.ini11}-\eqref{buk.ini13} obtained by Hybrid6 with exact solution \eqref{buk.ini53} at time $t=1$.}
 \label{buc.ex1}
\end{figure}

Again consider the same flux \eqref{buk.ini11} with the following initial data
\begin{equation}\label{buk.ini13}
u(x,0)=\left\{\begin{array}{ll}
                1 ~~~~   x >0.25, \\
                0 ~~~~  x\leq 0.25
                \end{array}\right.
\end{equation}
whose exact solution is given by
 \begin{equation}\label{buk.ini53}
u(x,t)=\left\{\begin{array}{ll}
                1 &~~~~  x<\dfrac{((\sqrt{3}-1)t+1)}{4},\medskip \\
                \dfrac{1}{2}\Big(\dfrac{t-4x+1}{t}\Big)& ~~~\dfrac{((\sqrt{3}-1)t+1)}{4}<x<\dfrac{t+1}{4},\medskip\\
                0&~~~x>\dfrac{t+1}{4}.
                \end{array}\right.
\end{equation}
Figure \ref{buc.ex1}(b), shows that Hybrid6 resolves the composite wave solution at time $t=1$ with non-oscillatory property.  Here, we have taken the CFL number to be 0.2 and the number of mesh points as 200.
\qed}
\end{example}
Finally, we test the Hybrid6 scheme with the one dimensional system of Euler equations of motion.
\begin{example}[1D System: Euler Equtions]{\rm

One-dimensional system of Euler equations can be written in the form of conservation laws \eqref{hyp.1d},
where $\bfm{u}$ and $\bfm{f}(\bfm{u})$ are vector of conservative variables and fluxes given respectively by\\
$$\bfm{u}=\begin{bmatrix}
    \rho \\
\rho u \\
e
   \end{bmatrix},
\hspace{1.0 cm}
\bfm{f}(\bfm{u})=\begin{bmatrix}
       \rho u\\
\rho u^2+p\\
e(u+p)
      \end{bmatrix}.$$
      Here, the dependent variables $\rho$, $p$, and $u$ denote density, pressure, and velocity, respectively. The energy $e$ is given by the relation
 \begin{equation}
e=\frac{p}{(\gamma-1)} +\frac{1}{2}\rho u^2,
\end{equation}
where $\gamma$ is the ratio of specific heat. 

We extend the Hybrid6 scheme for system of equations in component-wise manner. In order to make the algorithm more efficient, Dewar {\it et al.} \cite{dew-etal_15a} used only the WLTE of pressure term. This can be following in our case as well, but we use the WLTE for each equation separately.

The non-oscillatory nature of the Hybrid6 scheme and the resolution of
discontinuities are assessed on the one-dimensional Riemann problems given by
\begin{equation}
\bfm{u}(x,0)=\left\{\begin{array}{ll}\bfm{u}_l ~~~~~~~~x<x_0,\\
\bfm{u}_r ~~~~~~~~x>x_0.
\end{array}\right.
\end{equation}
These problems consist of two constant states $u_l$ and $u_r$ seperated with an initial discontinuity and involve the formation of a rarefaction wave, a
contact discontinuity, and a shock wave. For the numerical experiments, we consider the Sod's shock tube problem and the Lax's shock
tube problem. The exact solutions to these
problems are obtained using a Riemann solver (for more detail see Toro \cite{tor_09a}).
\begin{figure}
 \includegraphics[width=0.52\textwidth]{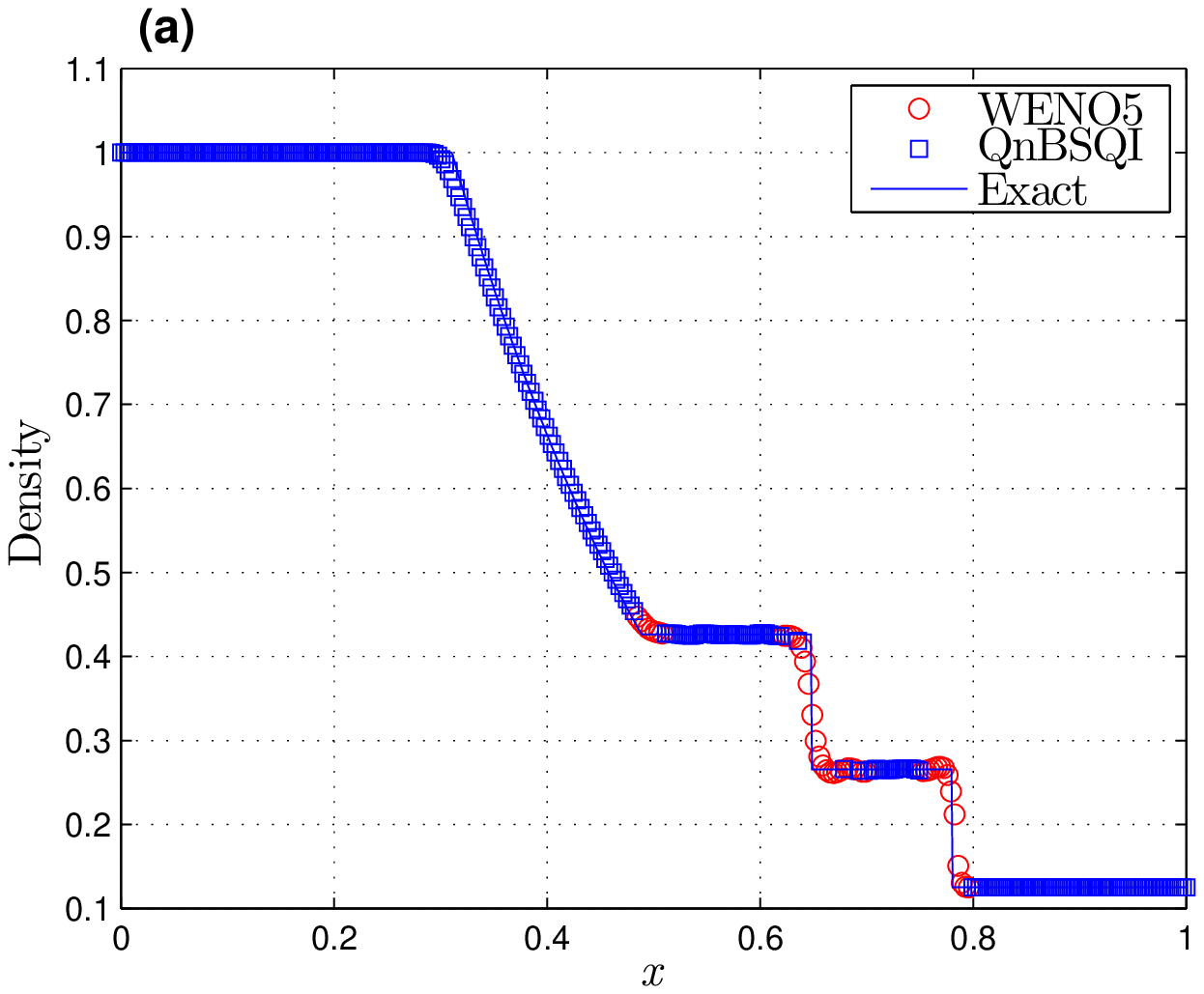}
 \includegraphics[width=0.52\textwidth]{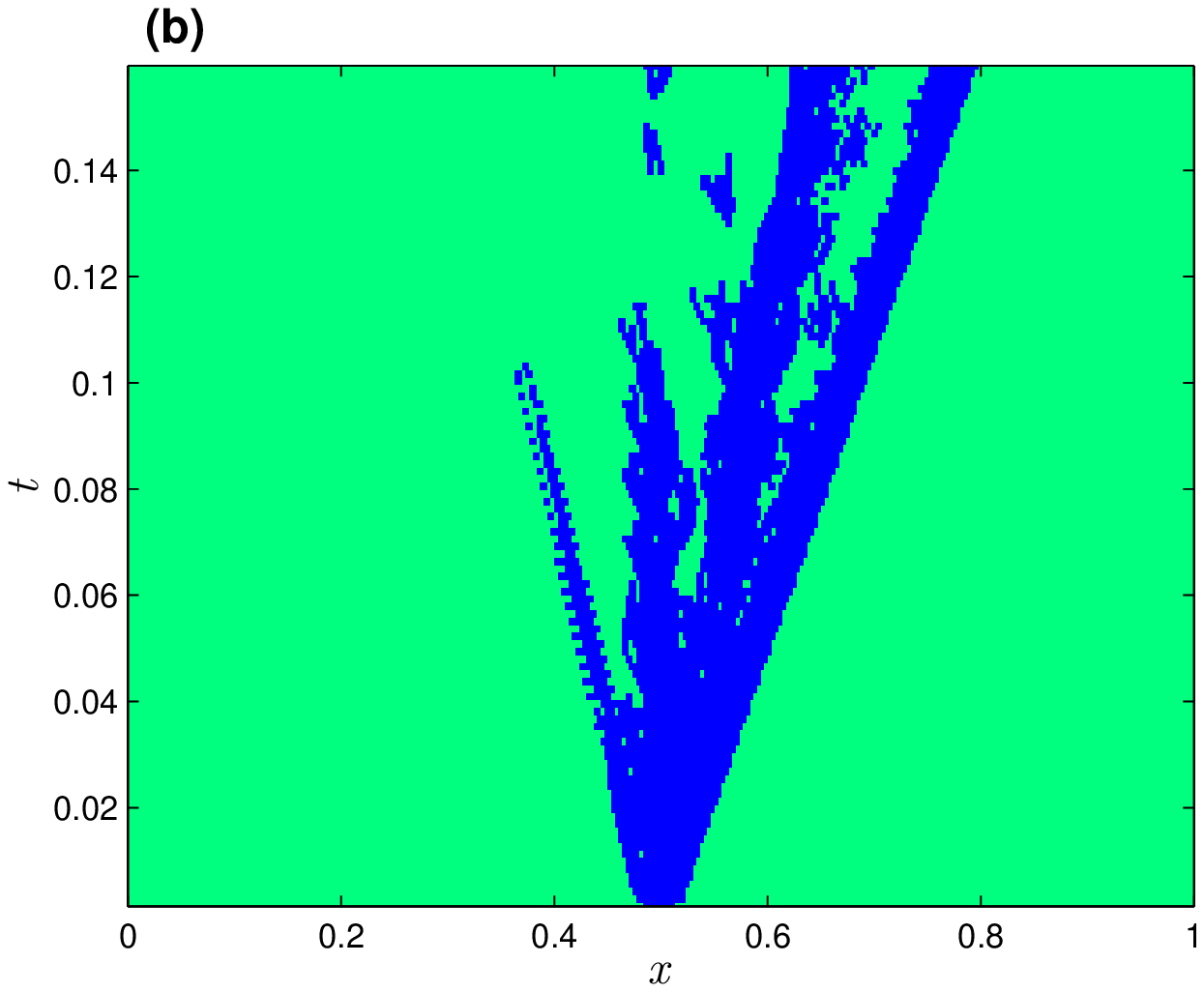}
  \includegraphics[width=0.52\textwidth]{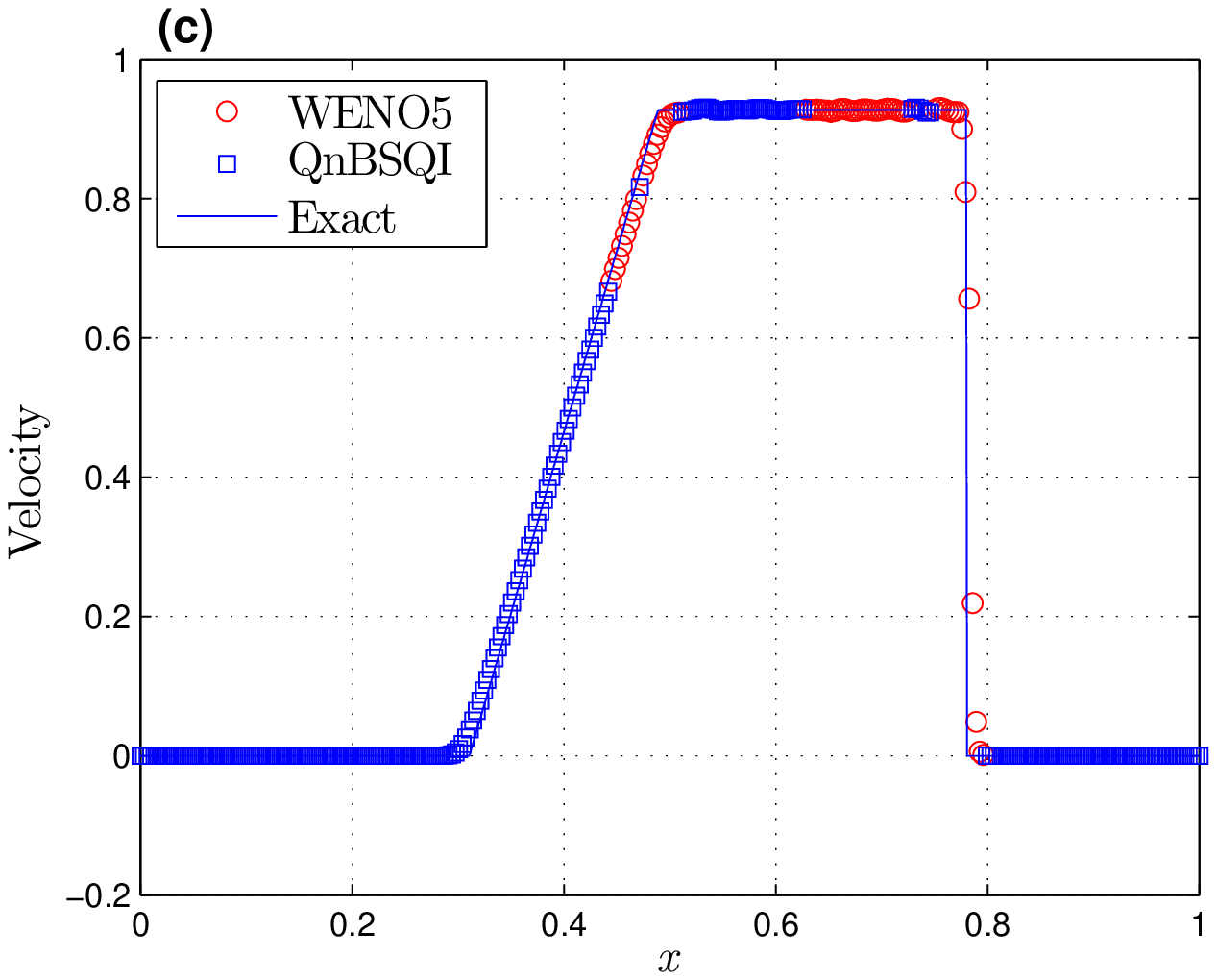}
 \includegraphics[width=0.52\textwidth]{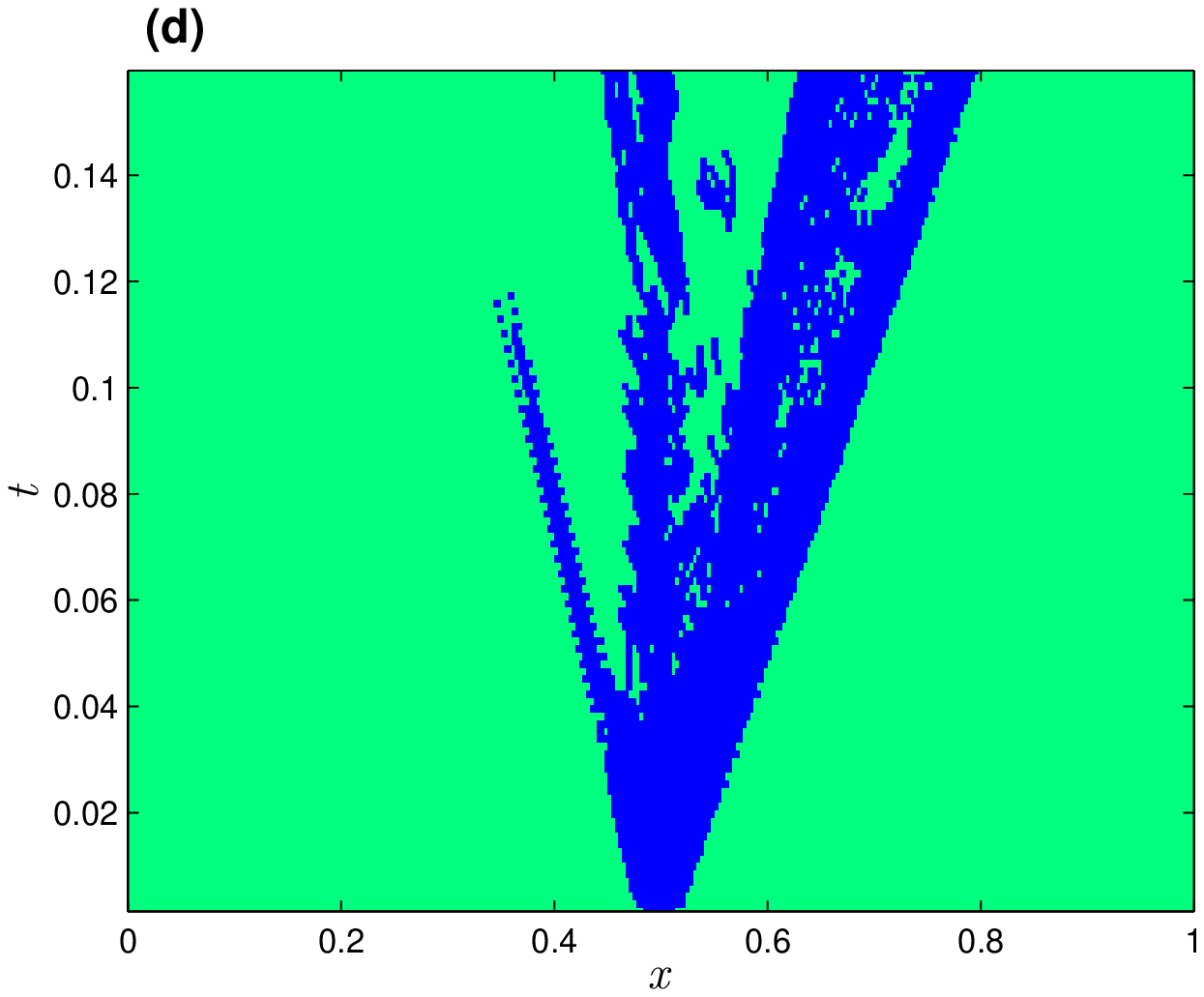}
 \includegraphics[width=0.52\textwidth]{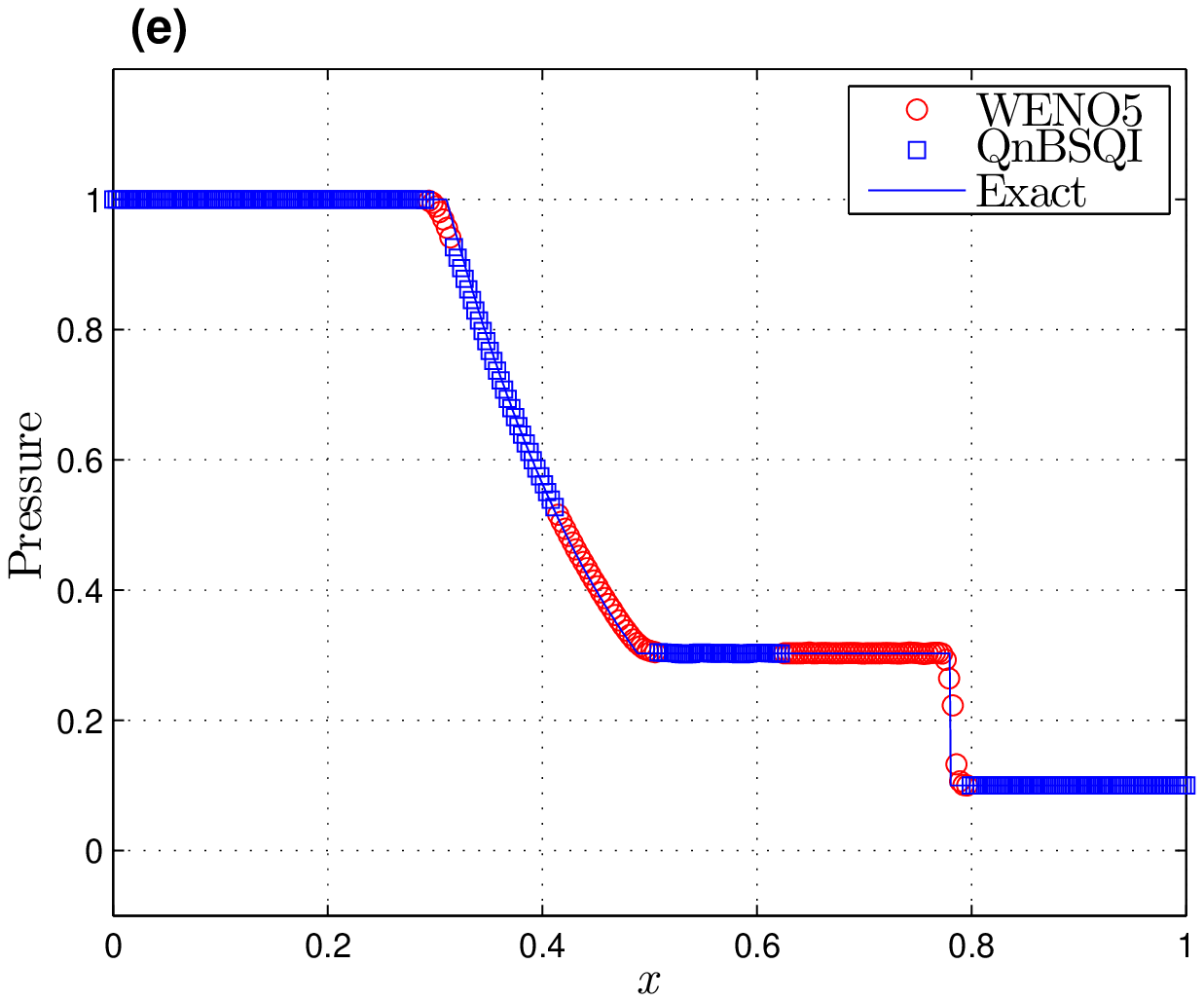}
 \includegraphics[width=0.52\textwidth]{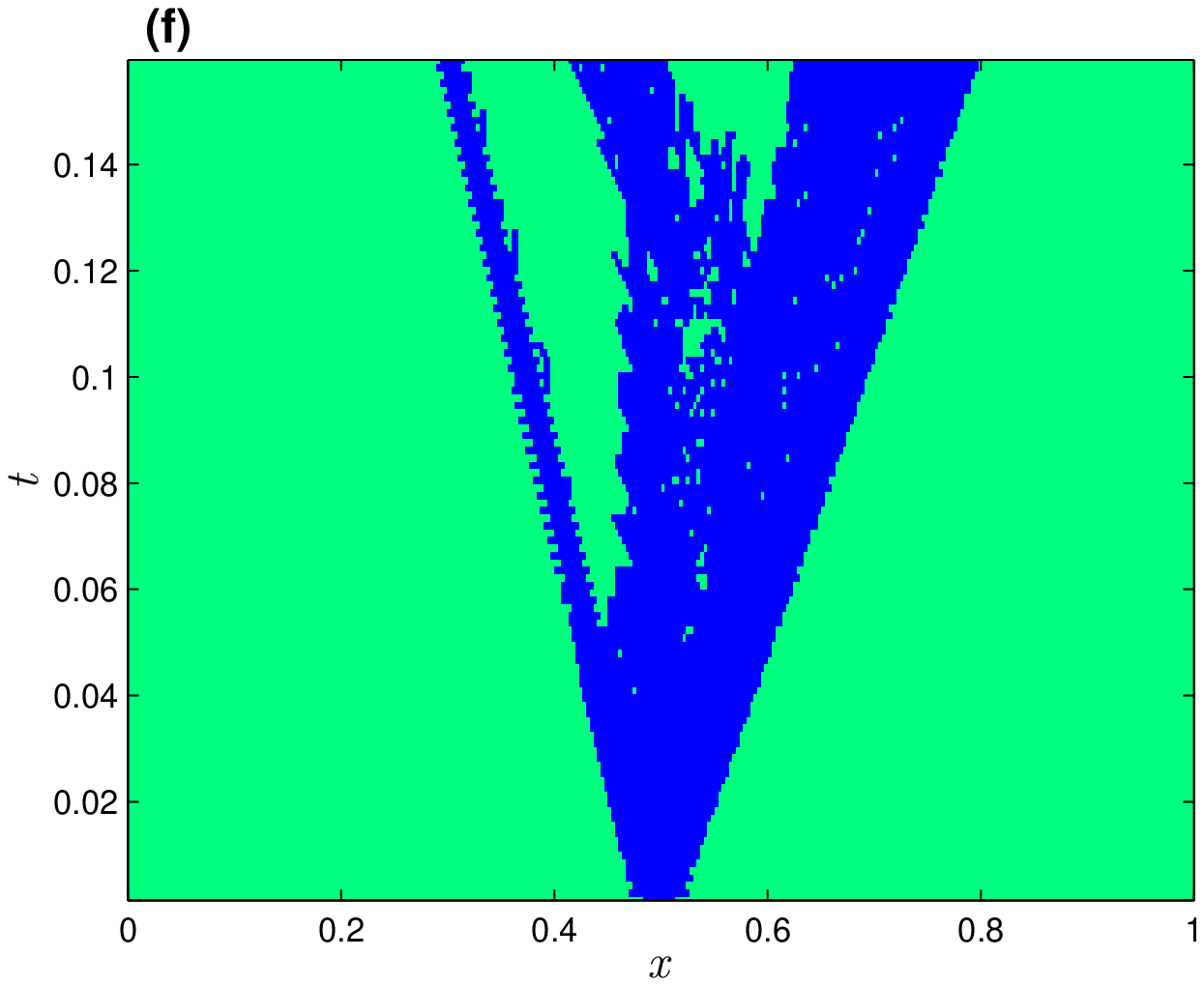}
 \caption{Comparison of the numerical solution at time $t=0.25$ obtained using Hybrid6 with the 
 reference solution, in the case of Sod's shock tube problem,  (a) density (c) velocity (e) pressure and corresponding plot of smooth indicator of the conservative variables for time interval $[0, 0.25]$  are shown in (b), (d), and (f), respectively. }
 \label{sod.34}
\end{figure}

Now we perform the numerical experiments to analyze Hybrid6 scheme for the system of Euler equations. We consider the
 Sod's shock tube problem arising from the following initial condition
\begin{equation}\label{Sod.system1}
(\rho, u, p)(x,0)=\left\{\begin{array}{ll}
                (1.0,0.0,1.0) ~~~~   x<0.5, \\
                (0.125,0.0,0.1) ~  x>0.5.
                \end{array}\right.
\end{equation}
Figure \ref{sod.34}(a), (c), and (e) depicts the comparison of the Hybrid6 results of the primitive variables density, velocity, and pressure along with the exact values.  The computation results are obtained over the domain $[0,1]$ at time $t=0.25$ with 300 mesh points. The CFL is kept constant equal to $0.3$. We observe from these results that the Hybrid6 scheme captures the shock with mild oscillation, this may be due to the transition of the scheme at shock position or WENO5 itself has mild oscillation at shock positions (see Kurganov and Liu \cite{kur-liu_12a}). The oscillation due to transition can be avoided by choosing a relatively smaller value of $K$ in the estimate \eqref{sind.est}. In Figure \ref{sod.34}(b), (d), and (f), we show the track of the smooth indicator with time. Initially smooth indicator starts with discontinuity position and detect the discontinuities and rarefaction tip along the time. 
\begin{figure}
 \includegraphics[width=0.52\textwidth]{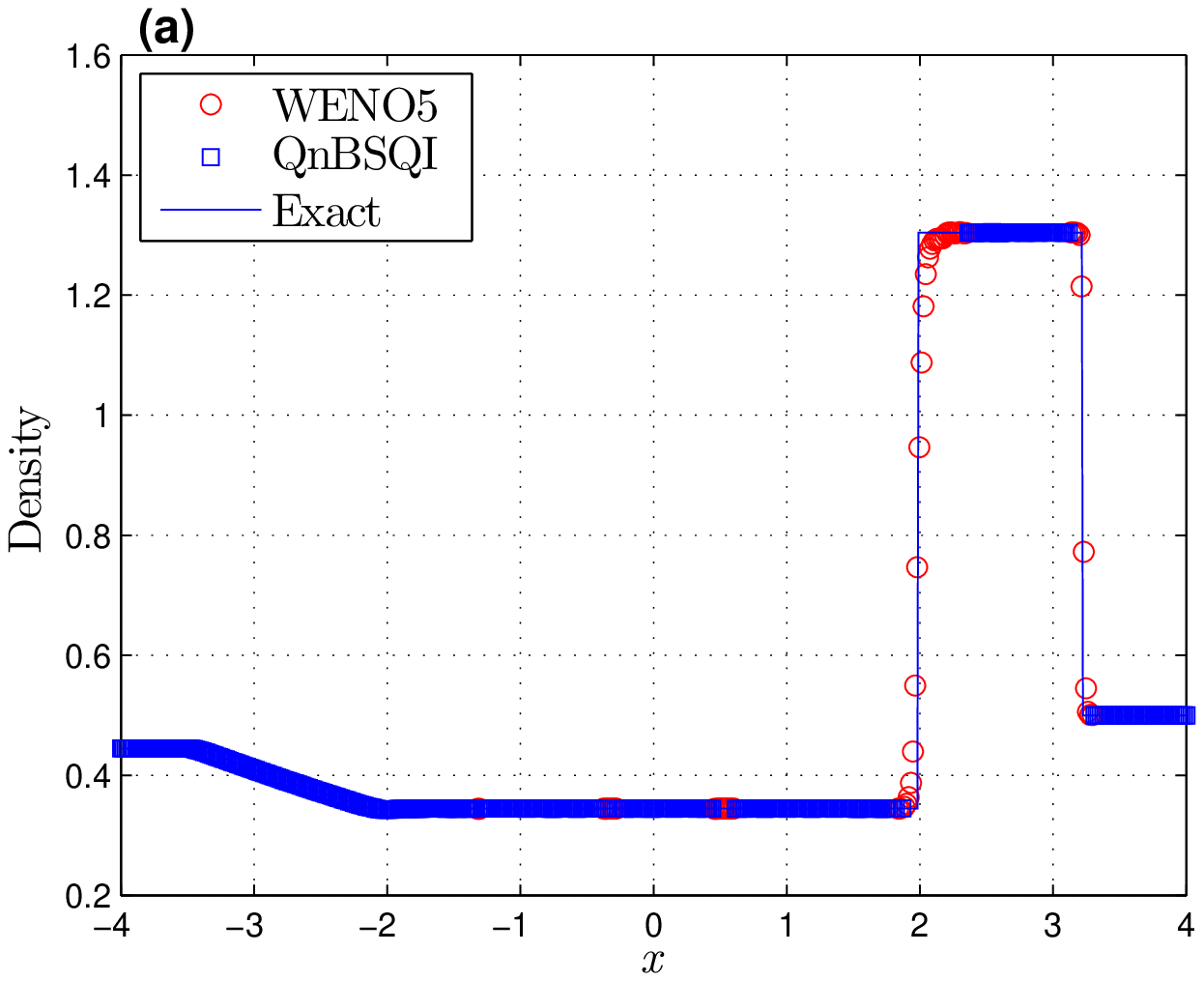}
 \includegraphics[width=0.52\textwidth]{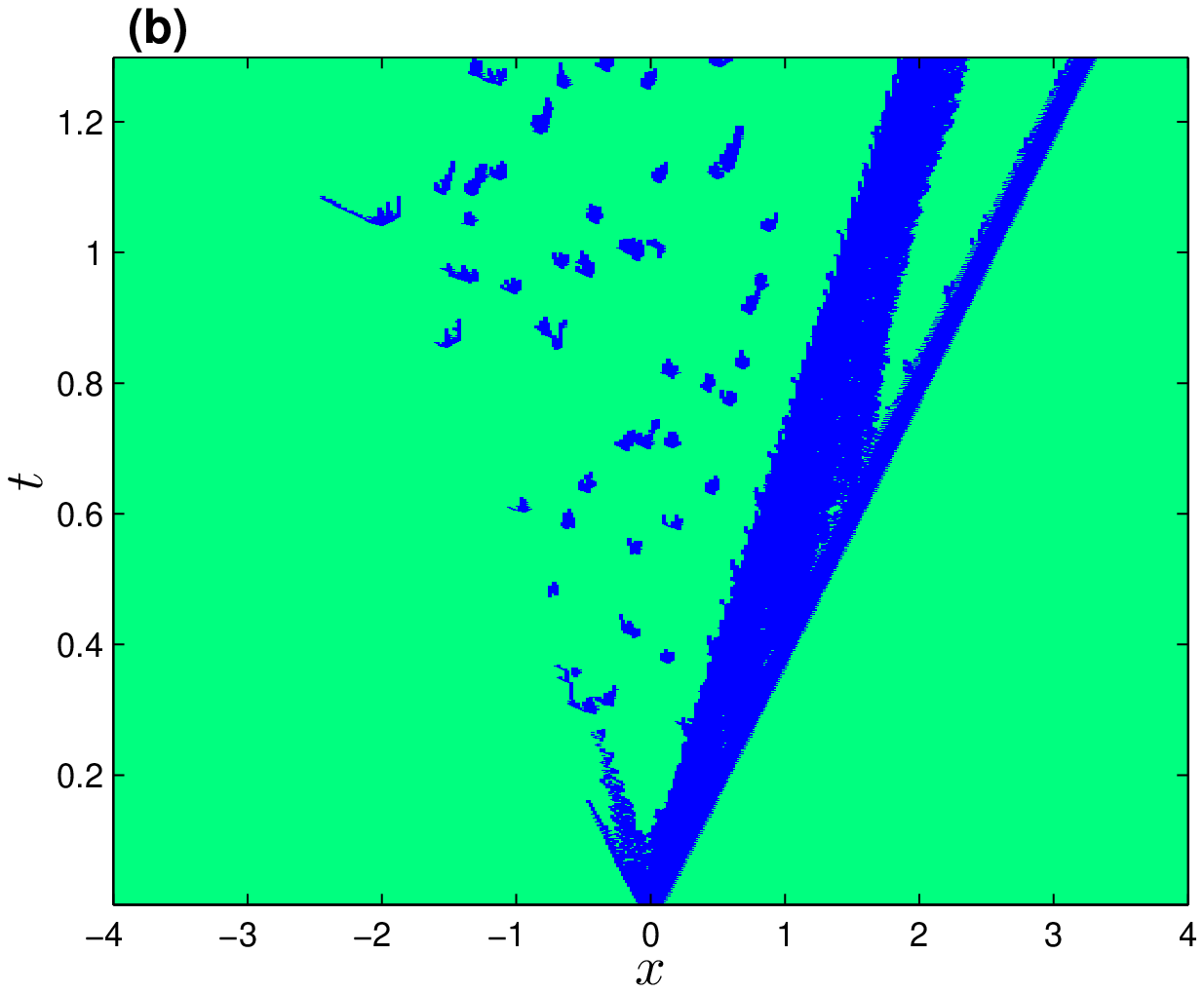}
 
 \caption{(a) Comparison of the numerical solution for density  at time $t=1.3$ obtained using Hybrid6 with the reference solution, in the case of Lax Problem. (b) Plot of smooth indicator for time interval $[0, 1.3]$.}
 \label{lax.34}
\end{figure}

 Further, we test Hybrid6 scheme with the Lax problem
   \begin{equation}\label{lax.system}
(\rho, u, p)(x,0)=\left\{\begin{array}{ll}
                (0.445,0.6980,3.528) ~~~~   x<0, \\
                (0.5,0.0,0.571) ~ ~~~~~~~~~~~ x>0.
                \end{array}\right.
\end{equation}
The numerical solution is computed over the domain $[-4,4]$ with 500 mesh points with the CFL number 0.4.
In Figure \ref{lax.34}, we have compared the density with exact solution at time $t=1.3$ and shown the plot of density smooth indicator over the $(x,t)$-plane.   
\qed}
\end{example}
\section{Order of Accuracy and Efficiency of the Hybrid Method}\label{ordeff.sec}

We have observed from the numerical experiments that the usage of the WENO
 scheme (in percentage) is negligible in the Hybrid6 scheme, especially when the number of grid points increases. However, the computational cost may increase due to the point-wise evolution of WLTE and the accuracy may decrease due to the shift from one method to other. These views demand further numerical tests for accuracy and efficiency of the developed hybrid schemes.  In this section, we extend our numerical experiments and perform some test cases to justify that the developed hybrid schemes are as accurate as the corresponding WENO schemes in the case when the solution is not smooth.  Also, we compare the performance of hybrid schemes with pure WENO schemes in term of CPU time taken to compute the solution.
\begin{table}[t]
\caption{Comparison of the WENO3 and  Hybrid4  in terms of $L^1$-error and their convergence rate in case of Burgers' equation with initial data  \eqref{buk.ini1} at time $t=0.5$.}  \label{table11}
\begin{center}
 \begin{tabular}{ |c | c | c |c|c| }
      \hline
 $N$ &      \multicolumn{2}{|c|}{WENO3}&\multicolumn{2}{|c|}{Hybrid4}\\ \cline{2-5}
 & $L^1$-error&  order &    $L^1$-error&  order  \\ \hline 
 50& 0.0670  &-  & 0.0662    & --  \\ \hline
 100&0.0294  & 1.18 & 0.0287     & 1.208  \\ \hline
 150&0.0174&   1.29&  0.0168  & 1.313 \\ \hline
 200&0.0132  & 0.97 &0.0127    & 0.982 \\ \hline

 \end{tabular}
\end{center}
\end{table}
\begin{table}[t]
\caption{Comparison of the WENO5 and  Hybrid6  in terms of $L^1$-error and their convergence rate in case of Burgers' equation with initial data  \eqref{buk.ini1} at time $t=0.5$.}  \label{table12}
\begin{center}
 \begin{tabular}{ |c | c | c |c|c| }
      \hline
 $N$ &      \multicolumn{2}{|c|}{WENO5}&\multicolumn{2}{|c|}{Hybrid6}\\ \cline{2-5}
 & $L^1$-error&  order &    $L^1$-error&  order  \\ \hline 
 50& 0.0620  &-  & 0.0610    & --  \\ \hline
 100&0.0267  & 1.215 & 0.0262     & 1.2186  \\ \hline
 150&0.0156&   1.326&  0.0152  & 1.3376 \\ \hline
 200&0.0118  & 0.969 &0.0116     & 0.9630 \\ \hline
 \end{tabular}
\end{center}
\end{table}

For the order of accuracy, we perform the numerical test in the case of Burgers equation with initial data \eqref{buk.ini1}.  The numerical solutions obtained using WENO scheme and the corresponding hybrid scheme are compared with the exact solution at time $t=0.5$ with CFL 0.1.  In Tables \ref{table11} and \ref{table12}, we show the $L^1$-error in WENO and the corresponding hybrid schemes are comparable and also the order of accuracy.  Recall that in Example \ref{Burger.ord.ex}, we have seen that when the solution is smooth, the CBSQI and QnBSQI schemes achieve the order of accuracy 4 and 6, respectively.  Since, the corresponding hybrid schemes use only the BSQI approximations as long as the solution remains smooth, we see that the developed hybrid schemes takes the maximum advantage of the accuracy of the BSQI schemes in the regions where the solution is smooth and are as accurate as the WENO schemes in the vicinity of the shocks.

 We now perform few test cases to study (numerically) the efficiency of the Hybrid4 and Hybrid6 schemes. 
For the efficiency test, we consider the Burgers' and the Buckley-Leverett equations with the initial condition \eqref{buk.ini1}, whereas in the case of Euler
 equations, we take the Sod shock tube test problem \eqref{Sod.system1}.  Tables \ref{efft.1}, \ref{efft.2}, and \ref{efft.3} list the CPU usage time in seconds of the WENO3, WENO5 schemes, and the Hybrid4 and Hybrid6 schemes for the above test cases, respectively.  These efficiency tests show that the Hybrid4 and Hybrid6 schemes whose numerical fluxes are given by \eqref{hybrid4.eq} and \eqref{hybrid6.eq}, respectively,  have the advantage of being considerably faster than the WENO schemes. 
We performed our computation on a computer with a i-5 processor and  a 4GB RAM.
 \begin{table}[t]
\centering
\caption{ Comparison of hybrid schemes for Burgers' equation with initial data \eqref{buk.ini1} with WENO3 and WENO5 at time $t=0.25$ 
in terms of CPU time taken (in seconds) with an increase in the number of mesh points.}
\begin{tabular}{|l| c|c| c| c|  c| r|}
\hline 
 $N$&Hybrid4&Hybrid6 & WENO3 & WENO5\\ 
\hline
100  & 0.001121      &0.001805          &0.002243    &0.004135\\
\hline
200  &  0.003848   &0.005410         &0.006815     &0.022249\\
\hline
400  &  0.011942   &0.0161005        &0.030349     &0.066770\\
\hline
800  &  0.043337  &0.059413         &0.117569     &0.245321\\
\hline
1600 & 0.161478  &0.221670         &0.505622     &1.123064\\
\hline
3200 &  0.627653   &0.868712         &2.096837     &4.726269\\
\hline
6400 &  2.480405  &3.336897         &7.616418     &16.414773\\
\hline
12800 & 9.786336    &13.11141         &28.71890     &60.144301\\
\hline
\end{tabular}
\label{efft.1}
\end{table}
\begin{table}
\centering
\caption{Comparison of hybrid schemes for Buckley-Leveratte equation with initial data \eqref{buk.ini1} with WENO3 and WENO5 at time $t=0.25$ 
in terms of CPU time taken (in seconds) with an increase in the number of mesh points.}
\begin{tabular}{|l| c| c| c|  c| c| r|}
\hline 
 $N$& Hybrid4& Hybrid6 & WENO3 & WENO5\\ 
\hline
100  & 0.002195	    &0.003030         &0.003225    &0.005836\\
\hline
200  &  0.009246   &0.008730         &0.009986     &0.025095\\
\hline
400  &  0.019462   &0.027194        &0.043534     &0.080060\\
\hline
800  &  0.069872   &0.100701         &0.155369     &0.332549\\
\hline
1600 & 0.283178    &0.369381         &0.721080     &1.1503515\\
\hline
3200 &  1.029097   &1.401521         &3.190346     &6.691973\\
\hline
6400 &  3.98037  &5.505693     &11.422844&23.609234\\
\hline
12800 & 15.7056   & 21.505811        &42.4329     &85.480766\\
\hline
\end{tabular}
\label{efft.2}
\end{table}
\begin{table}
\centering
\caption{Comparison of hybrid schemes for Euler's  equation of system in case of Sod shock tube problem with WENO3 and WENO5 at time $t=0.25$ 
in terms of CPU time taken (in seconds) with an increase in the number mesh points.}
\begin{tabular}{|l| c| c| c|  c| c| r|}
\hline 
 $N$& Hybrid4& Hybrid6& WENO3 & WENO5\\ 
\hline
100  & 0.003256    &0.004555         &0.05782  &0.0013244\\
\hline
200  &  0.011631   &0.020165         &0.028531     &0.054535\\
\hline
400  &  0.044551   &0.052351     &0.096789    &0.211751\\
\hline
800  &  0.146095    &0.204967         & 0.372008     &0.821623\\
\hline
1600 & 0.547058  &0.768439     &1.477857    &3.276202\\
\hline
3200 &  2.080392  &2.863554         &5.871739     &13.048635\\
\hline
6400 &  8.019332    &11.754494         &23.353153     &52.067757\\
\hline
12800 & 31.629588   &47.306144       &93.72000    &208.11320\\
\hline
\end{tabular}
\label{efft.3}
\end{table}


\section{Conclusion}\label{sec.conclusion}
We developed BSQI based numerical schemes, namely CBSQI and QnBSQI, for hyperbolic conservation laws using cubic and quintic BSQI techniques, respectively.  We showed numerically that they achieve 4$^{\rm th}$ and 6$^{\rm th}$ order of accuracy in space in the case of smooth solution.  However, as expected, these schemes develop spurious oscillations near shocks and therefore not suitable for numerical approximation.  To overcome this difficulty, we proposed two hybrid methods, namely, Hybrid4 and Hybrid6 where we conjugate the CBSQI and QnBSQI schemes with the third and the fifth order WENO schemes, respectively.  We used a weak local truncation based estimate to detect the high gradient or shock regions of the numerical solution.  We use this information to capture shocks using WENO scheme, whereas the BSQI based scheme is used in the smooth regions. For the time discretization, we consider a strong stability preserving Runge-Kutta method of order three. We have shown that the develop hybrid schemes are as accurate as WENO schemes and maintains non-oscillatory behaviour at the shock positions. Also, our numerical efficiency tests show that the developed hybrid schemes take only (approximately) one fifth of the time taken by  the WENO schemes.
\begin{acknowledgements}
The first author acknowledges the financial support of Council of Scientific and Industrial Research (CSIR), Government of India (09/087(0592)/2010-EMR-I). 
\end{acknowledgements}


\bibliography{bibfile}

\end{document}